\documentclass[11pt]{amsart}
\usepackage[utf8]{inputenc}
\usepackage[top=1in, left=1in, right=1in, bottom=1in]{geometry}
\usepackage{fancyhdr}
\usepackage{color,amsmath,amssymb,mathtools,graphicx}
\usepackage{amsfonts}
\usepackage{libertine}
\usepackage[libertine]{newtxmath}
\usepackage{epsfig,fancyhdr}
\usepackage{dsfont} 
\usepackage{graphicx,float}
\usepackage{graphicx}
\usepackage{epsfig,color,fancyhdr,setspace}
\usepackage{epstopdf}
\usepackage{wrapfig}
\usepackage{soul}
\usepackage{import}
\usepackage{lineno}
\usepackage{stackrel}
\usepackage[dvipsnames]{xcolor}

\usepackage[hidelinks,pagebackref]{hyperref}

\usepackage[abs]{overpic}
\usepackage[center]{caption}
\usepackage{subcaption}
\usepackage{pst-node}
\usepackage{tikz-cd}
\usepackage{enumitem}
\counterwithin{figure}{section}
\usepackage{tikz-cd}
\usepackage{marvosym}
\usepackage{extarrows}
\usepackage{amsmath}
\usepackage{graphicx,amsmath,mathtools,float}
\usepackage{setspace}
\graphicspath{{./images/}}
\usepackage{stackrel}
\usepackage{yfonts}
\usepackage{enumitem, hyperref}\makeatletter
\def\namedlabel#1#2{\begingroup
    #2%
    \def\@currentlabel{#2}%
    \phantomsection\label{#1}\endgroup
}
\makeatother
\usepackage{type1cm} 
\usepackage{blindtext}
\usepackage{scalerel,stackengine}
\usepackage{tikz}
\usepackage{mathtools}
\usepackage{amssymb}
\usepackage{tikz}
\usetikzlibrary{shapes,decorations.pathmorphing}
\newcommand{\KP}[1]{%
	\begin{tikzpicture}[baseline=-\dimexpr\fontdimen22\textfont2\relax]
		#1
	\end{tikzpicture}%
}

\newcommand{\KPB}{%
	\KP{
		\draw[color=gray,thick] (-0.3,0.3) -- (0.3,-0.3);
		\draw[color=gray,thick] (-0.3,-0.3) -- (-0.05,-0.05);
		\draw[color=gray,thick] (0.05,0.05) -- (0.3,0.3);
	}%
}
\newcommand{\KPC}{%
	\KP{%
		\draw[color=gray,thick] (-0.3,0.3) .. controls (0,-0.05) .. (0.3,0.3);
		\draw[color=gray,thick] (-0.3,-0.3) .. controls (0,0.05) .. (0.3,-0.3);
	}%
}
\newcommand{\KPD}{%
	\KP{%
		\draw[color=gray,thick] (-0.3,-0.3) .. controls (0.05,0) .. (-0.3,0.3);
		\draw[color=gray,thick] (0.3,-0.3) .. controls (-0.05,0) .. (0.3,0.3);
	}%
}
\newcommand\quotient[2]{
	\mathchoice
	{
		\text{\raise1ex\hbox{$#1$}\Big/\lower1ex\hbox{$#2$}}%
	}
	{
		#1\,/\,#2
	}
	{
		#1\,/\,#2
	}
	{
		#1\,/\,#2
	}
}

\usepackage{tocbasic}
\newtheorem{theorem}{Theorem}[section]
\newtheorem{lemma}[theorem]{Lemma}
\newtheorem{definition}[theorem]{Definition}
\newtheorem{example}[theorem]{Example}

\title{ALEXANDER-CONWAY AND BRACKET POLYNOMIALS OF PRETZEL LINKS $\boldsymbol{P(1,1,n)}$}

\author{Alan Hern\'andez-Flores}
\address{Universidad Nacional Autónoma de Honduras UNAH}
\email{alanflores@unah.hn}

\author{Gabriel Montoya-Vega}
\address{Department of Mathematics, The Graduate Center CUNY, NY, USA, and \newline \indent Department of Mathematics, University of Puerto Rico-R\'io Piedras, San Juan, PR}
\email{gabrielmontoyavega@gmail.com}

\subjclass[2020]{Primary: 57K10. Secondary: 57K14.}
\keywords{Knots and links, pretzel links, Alexander-Conway polynomial, Kauffman bracket polynomial.}

\begin{document}
\begin{abstract}
Polynomial invariants constitute a dynamic and essential area of study in the mathematical theory of knots. From the pioneer Alexander polynomial, the revolutionary Jones polynomial, to the collectively discovered HOMFLYPT polynomial (just to mention a few), these algebraic expressions have been central to the understanding of knots and links. The introduction of Khovanov homology has sparked significant interest in the categorification of these polynomials, offering deeper insights into their topological and algebraic properties. In this work, we revisit two prominent polynomial invariants—the Alexander-Conway and the Kauffman bracket polynomials—and focus specifically on the polynomials associated with the family of three-strand pretzel links $P(1,1,n)$.  
\end{abstract}

\maketitle

\tableofcontents

\section{Introduction}
Modern knot theory traces its origins to a 1679 speculation by Leibniz, who proposed that alongside calculus and analytical geometry, there should be a geometry of position (geometria situs) focused solely on relationships determined by position, without regard to magnitudes. The first compelling example of geometria situs was provided by Leonard Euler in 1735, when he solved the famous problem of the bridges of Königsberg \cite{PBIMW}. Classical knot theory studies the embeddings of a circle (knots), or several circles (links), up to natural deformations in $\mathbb{R}^3$ highlighting the classification as a fundamental problem. This classification, or enumeration, is done up to the natural movement in space which is called an ambient isotopy. In 1927, Reidemeister showed that two link diagrams (possibly oriented) are isotopic if and only if they are connected by a finite sequence of moves, called Reidemeister moves, and planar isotopy; see \cite{Ada} for a detailed and friendly introduction to the mathematical theory of knots. Expanding on this concept, in order to distinguish knots and links we look for \textit{invariants of links} $-$ properties of links that remain unchanged under ambient isotopy. In other words, to show that a property $P(D)$ of a diagram $D$ is an invariant, one must confirm its persistence through Reidemeister moves. For instance, the Jones polynomial is an invariant of links \cite{Jon}. The article is organized in the following way. In Section
\ref{PretzelIntro} the standard presentation of pretzel links is recalled; in Section \ref{ACPoly} the Alexander-Conway polynomial is defined and a formula for its calculation is presented. Analogously, in Section \ref{KBPoly} the Kauffman bracket polynomial is introduced and a calculation formula is presented. Finally in Section \ref{Conclusions} we summarize the results and speculate about some possible future research ideas.

\section{Pretzel Links}\label{PretzelIntro}
\begin{figure}[H]
    \centering
\includegraphics[width=0.5\textwidth]{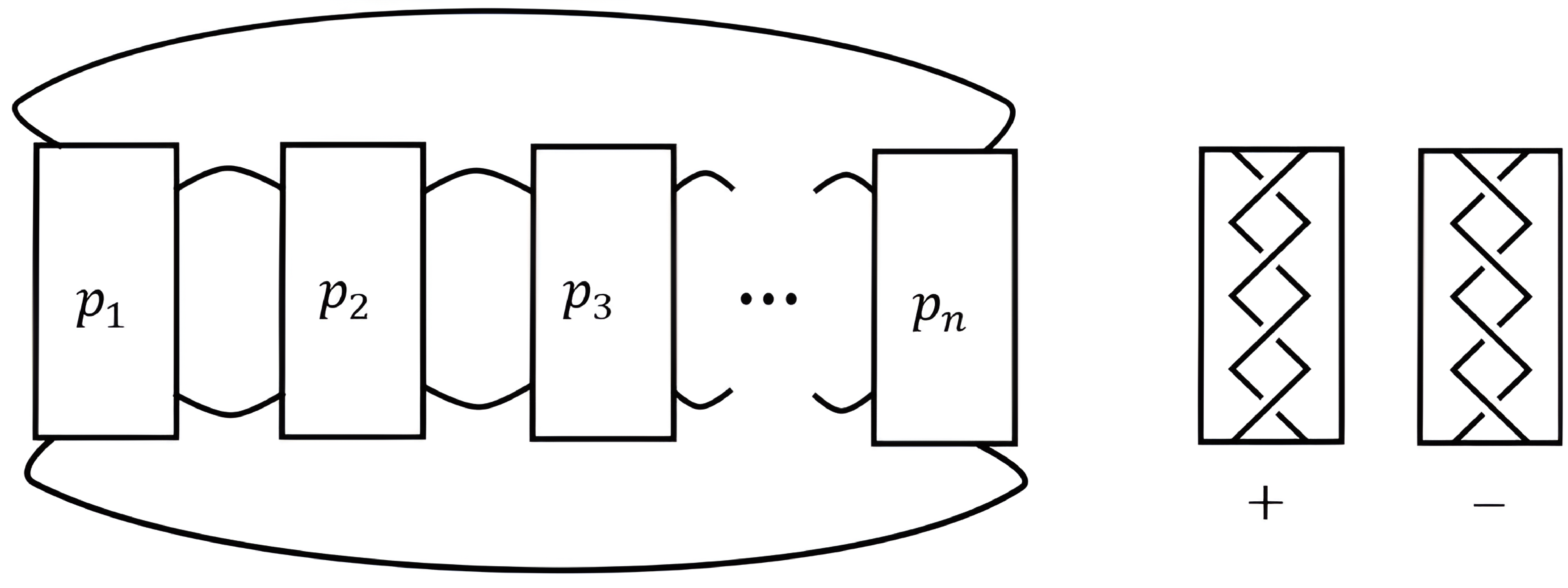}
    \caption{Pretzel link $P(p_1,p_2,p_3, \dots, p_n)$.}
    \label{DefiPretzel}
\end{figure}

Pretzel links constitute a well-studied family in the mathematical theory of knots. Recall that a pretzel link $P(p_1,p_2,p_3, \dots, p_n)$ is defined by an $n$-tuple of integers $(p_1,p_2,p_3, \dots ,p_n)$ so that each $p_i$, $i=1, \dots, n$, is different from zero. The number $p_i$ is the number of crossings in the tangle of the column $i$; see Figure \ref{DefiPretzel} showing the standard form of a pretzel link and the convention for positive and negative tangles. 
 
\section{Alexander-Conway Polynomial}\label{ACPoly}

In 1912, George David Birkhoff introduced the chromatic polynomial of a graph while attempting to prove the four-color problem, first posed by Francis Guthrie in 1852 \cite{Bir}. In simple terms, this polynomial counts the number of ways to label or color the vertices of a graph such that adjacent vertices receive different colors. A significant advancement in the theory of knot invariants came with the discovery of a Laurent polynomial invariant by James Waddell Alexander. Alexander, born in New Jersey, earned his PhD in 1915 from Princeton University under the mentorship of Oswald Veblen. His career overlapped with Birkhoff’s time at Princeton, suggesting that Alexander may have been familiar with Birkhoff’s chromatic polynomial. Alexander developed multiple approaches to this polynomial, including a combinatorial one and another involving the knot group. For a more detailed history of the polynomial, see, for instance, Chapter 2 of \cite{PBIMW}. In the 1960s, John Conway rediscovered Alexander's formula and introduced a normalized version of the polynomial, as outlined below.

\begin{figure}[H]
\centering
\begin{subfigure}{.3\textwidth}
\centering		\includegraphics[width=0.25\linewidth]{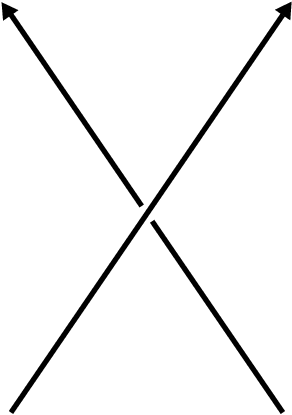}
\caption{$L_+$}
\end{subfigure}
\begin{subfigure}{.3\textwidth}
\centering
\includegraphics[width=0.25\linewidth]{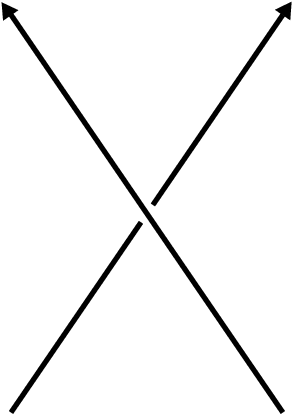}
\caption{$L_-$}
\end{subfigure}
\begin{subfigure}{.3\textwidth}
\centering		\includegraphics[width=0.25\linewidth]{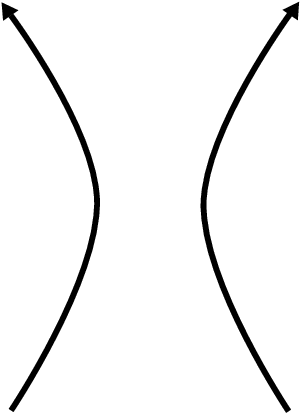}
\caption{$L_0$}
\end{subfigure}
\caption{Skein diagrams.}
\label{Skeinrelation}
\end{figure}

\begin{definition}[Alexander-Conway polynomial] \label{ACdefi}
    Given an oriented knot (or link) diagram $K$, we may assign to it a Laurent polynomial, $\nabla_K(z)$, by means of the following two axioms:

    \begin{itemize}
    \item[(i)] \(\nabla_{\bigcirc}(z) = 1\), where $\bigcirc$ is any diagram of the trivial knot.  
    \item[(ii)] Suppose \(L_{+}\), \(L_{-}\), and \(L_{0}\) are regular link diagrams that are identical except at a single crossing. At this crossing, the diagrams differ in the manner illustrated in Figure \ref{Skeinrelation}. The polynomials of these three diagrams are related by the following equation:  
    \begin{gather*}
        \nabla_{L_{+}}(z) = \nabla_{L_{-}}(z) + z \nabla_{L_{0}}(z).
    \end{gather*}
\end{itemize}
\end{definition}

\begin{example}
We calculate the Alexander-Conway polynomial of the right-trefoil knot. First we take a crossing of the knot and determine the  diagrams $L_+, L_-$ and $L_0$. Observe that, $L_+$ is the trefoil, $L_-$ is equivalent to the unknot, and $L_0$ is the Hopf link; see Figure \ref{SkeinrelationTrefoilExample}.

\begin{figure}[H]
\centering
\begin{subfigure}{.3\textwidth}
\centering		\includegraphics[width=0.6\linewidth]{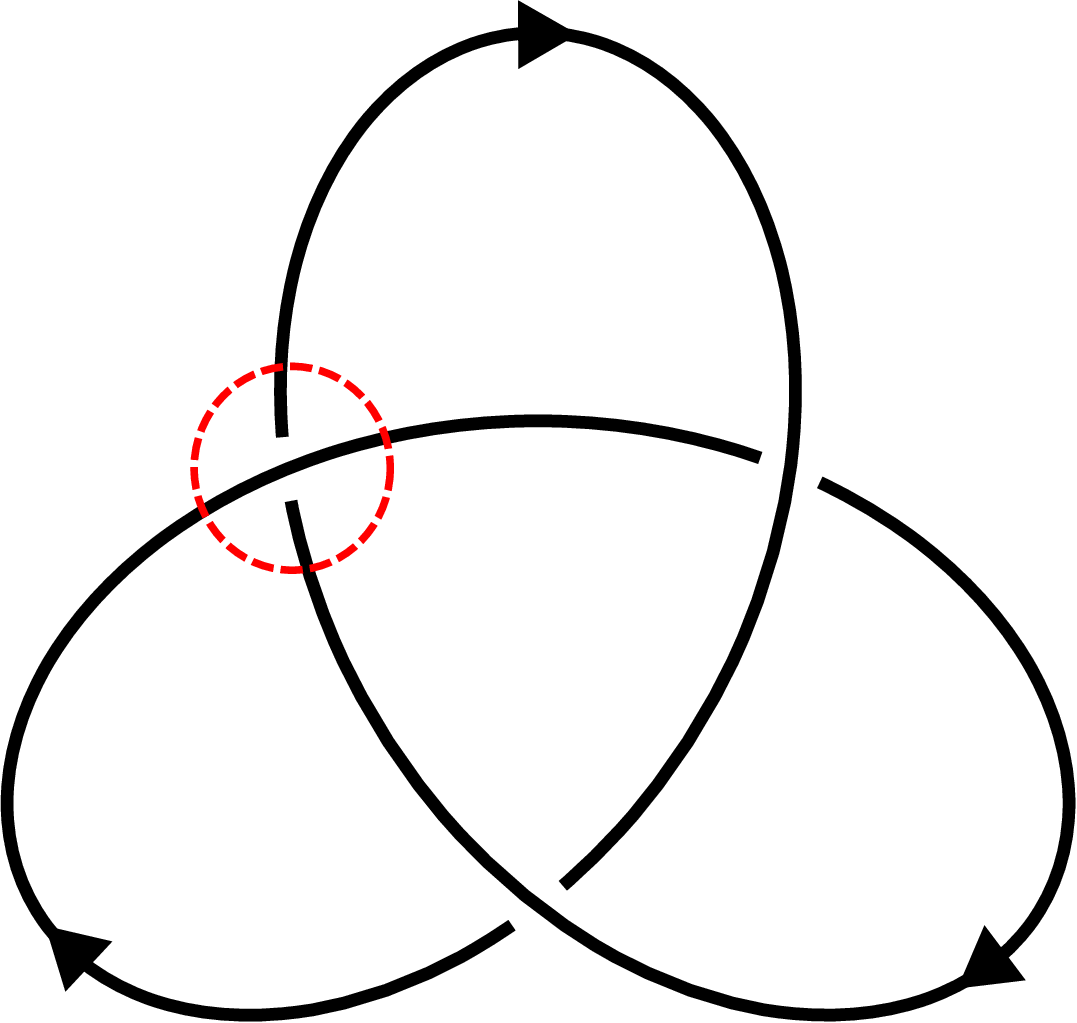}
\caption{$L_+$}
\end{subfigure}
\begin{subfigure}{.3\textwidth}
\centering
\includegraphics[width=0.6\linewidth]{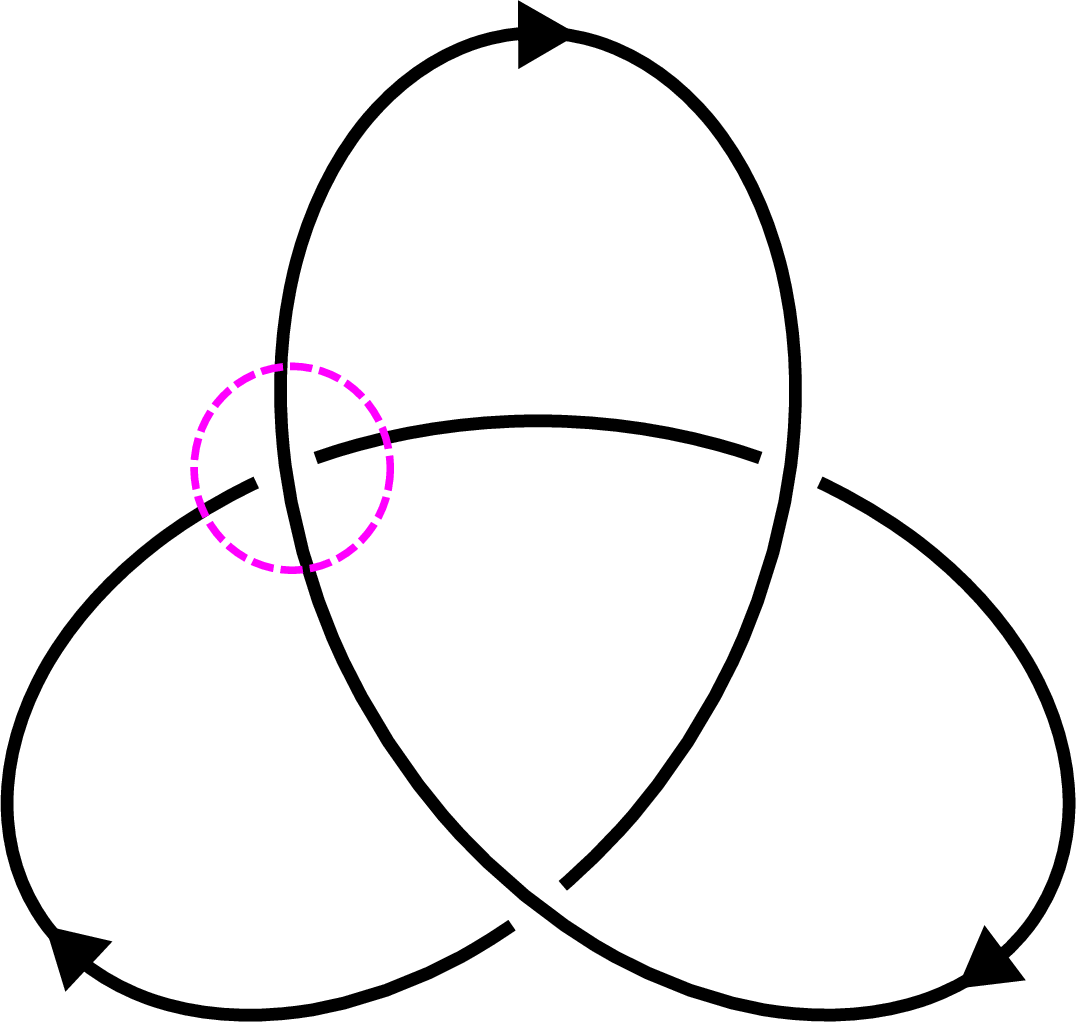}
\caption{$L_-$}
\end{subfigure}
\begin{subfigure}{.3\textwidth}
\centering		\includegraphics[width=0.6\linewidth]{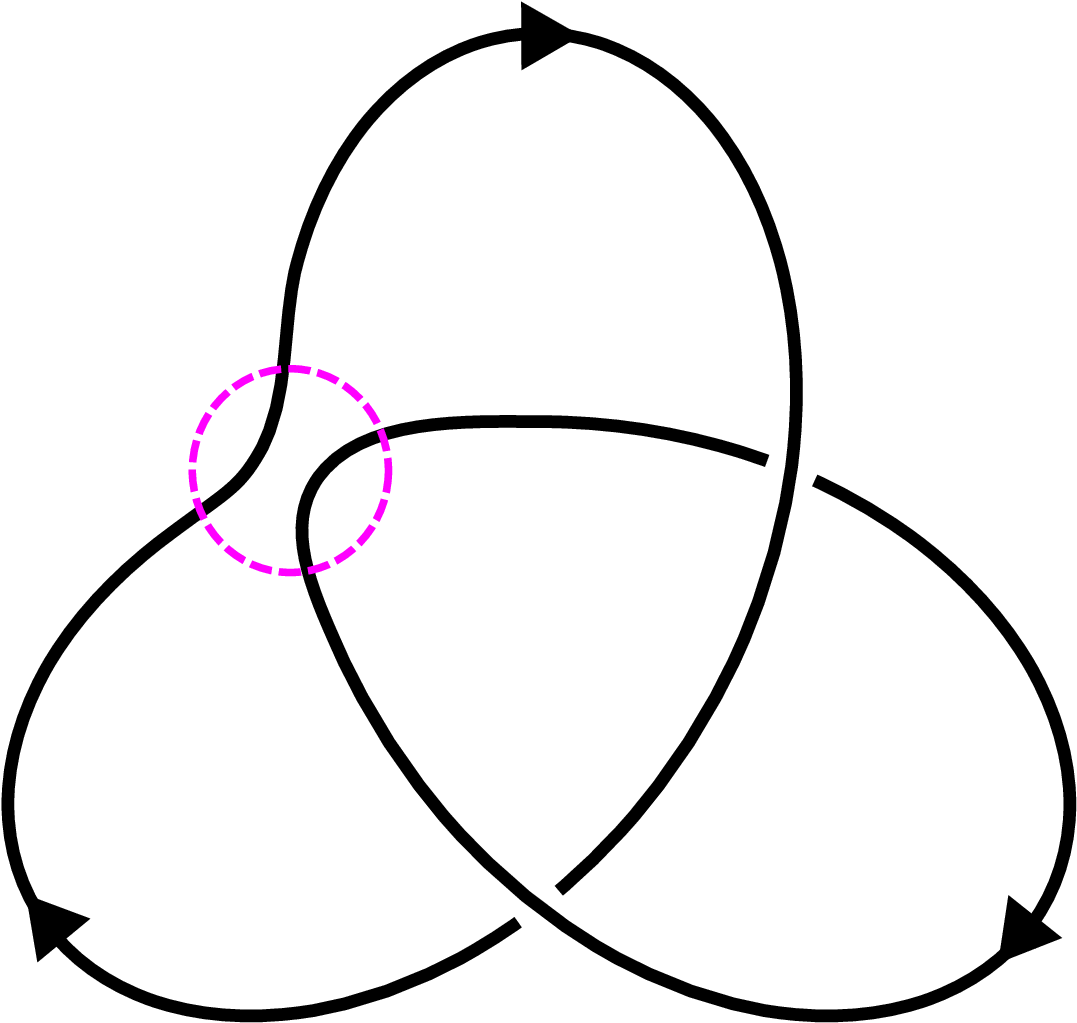}
\caption{$L_0$}
\end{subfigure}
\caption{Trefoil knot skein diagrams.}
\label{SkeinrelationTrefoilExample}
\end{figure}

Now, we do the same process for $L_0$ which is the Hopf link in Figure \ref{SkeinrelationHopfExample}, where $\widehat{L_+}$ is the Hopf link, $\widehat{L_-}$ is the trivial link of two components, and $\widehat{L_0}$ is equivalent to the unknot.

 \begin{figure}[H]
\centering
\begin{subfigure}{.3\textwidth}
\centering		\includegraphics[width=0.6\linewidth]{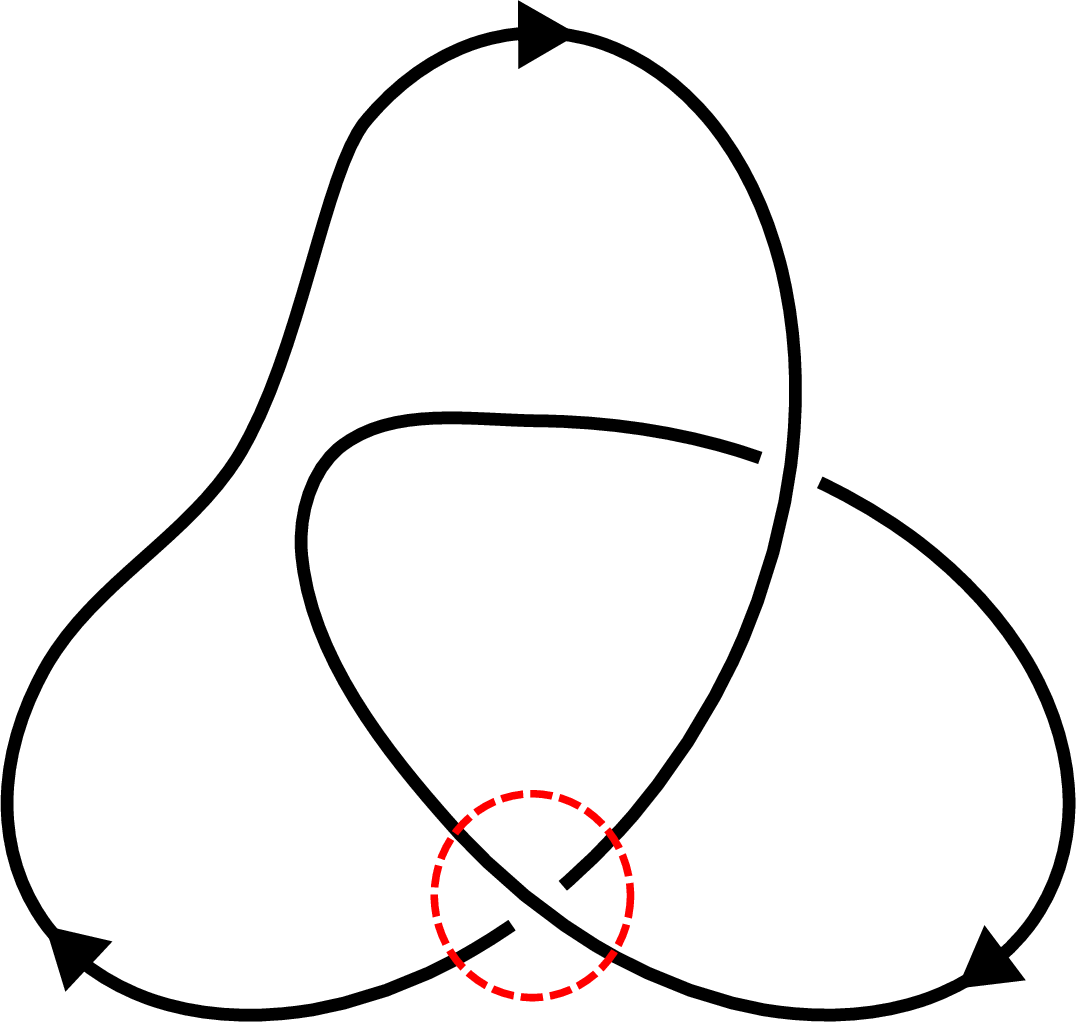}
\caption{$\widehat{L_+}$}
\end{subfigure}
\begin{subfigure}{.3\textwidth}
\centering
\includegraphics[width=0.6\linewidth]{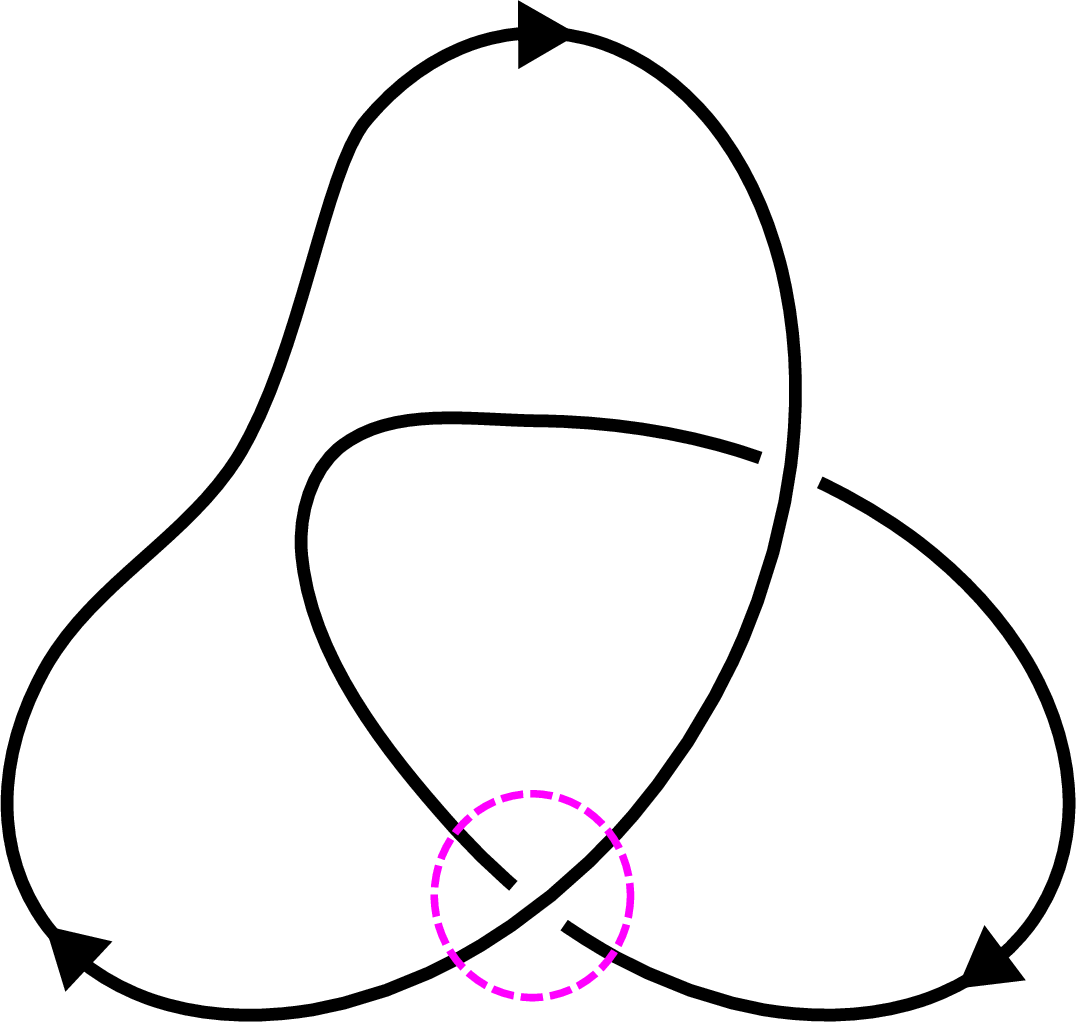}
\caption{$\widehat{L_-}$}
\end{subfigure}
\begin{subfigure}{.3\textwidth}
\centering		\includegraphics[width=0.6\linewidth]{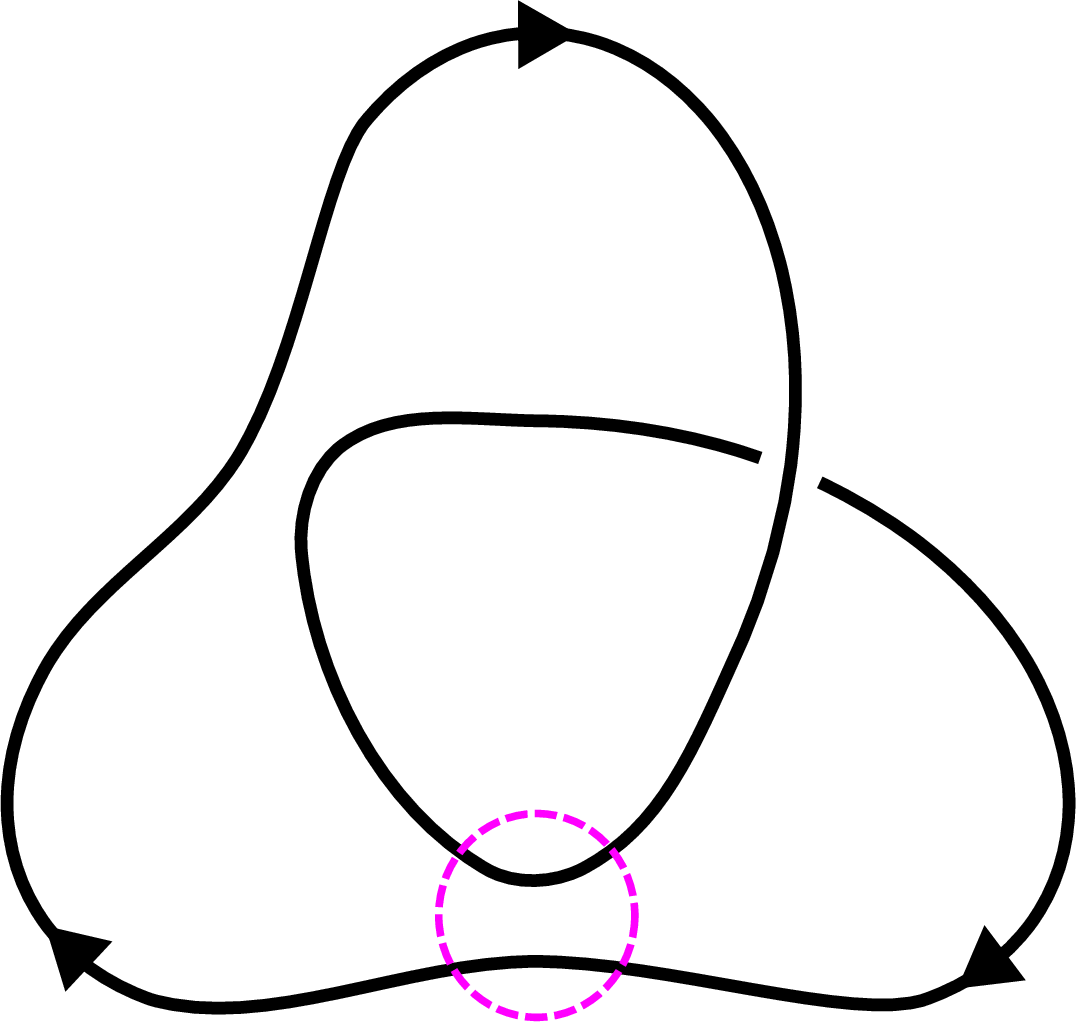}
\caption{$\widehat{L_0}$}
\end{subfigure}
\caption{Hopf link skein diagrams.}
\label{SkeinrelationHopfExample}
\end{figure}

Thus, the Alexander-Conway polynomial of the right-trefoil knot is given by:
\begin{gather*}
    \nabla_{L_+}(z)=\nabla_{L_-}(z)+z\nabla_{L_0}(z)\\
    \hspace{1.2cm}=\nabla_{\bigcirc}(z)+z\nabla_{\widehat{L_+}}(z)\\
    \hspace{3.1cm}=\nabla_{\bigcirc}(z)+z[\nabla_{\widehat{L_-}}(z)+z\nabla_{\widehat{L_0}}(z)]\\
    \hspace{3.2cm}=\nabla_{\bigcirc}(z)+z[\nabla_{\bigcirc \bigcirc}(z)+z\nabla_{\bigcirc}(z)]\\
    =1+z[0+z]\\
    \nabla(z)=1+z^2.
\end{gather*}
In the previous sequence of equalities we used the fact that $\nabla_{\bigcirc \bigcirc}(z)=0$.
\end{example}

\subsection{Alexander-Conway Polynomial of Pretzel links \texorpdfstring{$\boldsymbol{P(1,1,n)}$}{P(1,1,n)}.}
In this section, we present a formula for the calculation of the Alexander-Conway polynomial of pretzel links of the form \( P(1,1,n) \), where \( n \) is an integer.

\begin{theorem}
    The Alexander-Conway polynomial of a pretzel link \( P(1,1,n) \) with \( n \in \mathbb{Z} \) is given by:
\begin{gather}
    \nabla_{P(1,1,n)}(z) = \left\{ \begin{array}{lcc} 1 + \left( \frac{n+1}{2} \right) z^2, & \text{if} & n \hspace{0.1cm} \text{is odd} \\ \\ 1 - \left( \frac{n}{2} \right) z^2, & \text{if} & n \hspace{0.1cm} \text{is even.} \\ \end{array} \right.
\end{gather}
\label{the_AC}
\end{theorem}

\begin{proof}
Suppose that \( n \) is a positive odd integer, then we can write \( n = 2q + 1 \), where \( q \in \mathbb{Z}^+ \cup \{0\} \). Let us apply induction on \( q \). For \( q = 0 \), we have \( n = 1 \), which corresponds to the knot \( P(1,1,1) \), and it is equivalent to the trefoil knot. Therefore, we have
\[
    \nabla_{P(1,1,1)}(z) = 1 + z^2 = 1 + \left( \frac{1+1}{2} \right) z^2.
\]
Assume the statement is true for \( q = k \), i.e., for \( P(1,1,2k+1) \). Consider \( q = k + 1 \). In this case, \( n = 2(k+1) + 1 \), and we will prove the formula for \( P(1,1,2(k+1)+1) \). Let us take a crossing of the knot \( P(1,1,2(k+1)+1) \); we will select the first crossing of the third tangle and now, we determine the diagrams \( L_+ \), \( L_- \), and \( L_0 \). Notice that \( L_- \) is the knot \( P(1,1,2(k+1)+1) \). In \( L_+ \), due to the crossing change, two half twists are undone in the third tangle, which means that \( L_+ \) is equivalent to the knot \( P(1,1,2k+1) \), and by the induction hypothesis, we have \( \nabla_{P(1,1,2k+1)} = 1 + \left( \frac{(2k+1)+1}{2} \right) z^2 \). Finally, \( L_0 \) is the left-handed Hopf link (and recall that $\nabla_{L_0}(z)=-z$); see Figure \ref{Skeinrelationproof1_1}. Using the skein relation, we find that the polynomial is:

\begin{equation*}
    \begin{split}
        \nabla_{P(1,1,n)}(z)=\nabla_{P(1,1,2(k+1)+1)}(z)&=\nabla_{L_-}(z)\\&=\nabla_{L_+}(z)-z\nabla_{L_0}(z)\\
        &=\nabla_{P(1,1,2k+1)}(z)-z\nabla_{L_0}(z)\\
        &=1+\left(\frac{(2k+1)+1}{2}\right)z^2-z(-z)\\
        &=1+\left(\frac{(2k+1)+1}{2}\right)z^2+z^2\\
        &=1+\left(\frac{(2k+1)+1}{2}+1\right)z^2\\
        &=1+\left(\frac{(2(k+1)+1)+1}{2}\right)z^2\\
        &=1+\left(\frac{n+1}{2}\right)z^2.\\
    \end{split}
\end{equation*}

\begin{figure}[ht]
\centering
\begin{subfigure}{.3\textwidth}
\centering		\includegraphics[width=0.7\linewidth]{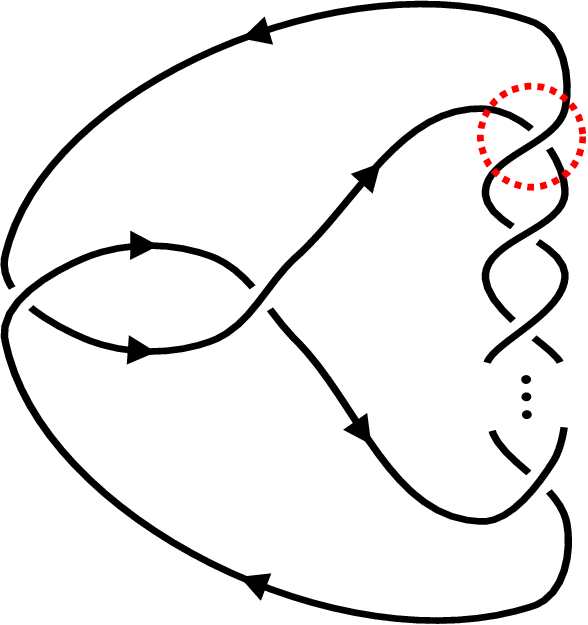}
\caption{$L_-$}
\label{Lminus}
\end{subfigure}
\begin{subfigure}{.3\textwidth}
\centering
\includegraphics[width=0.7\linewidth]{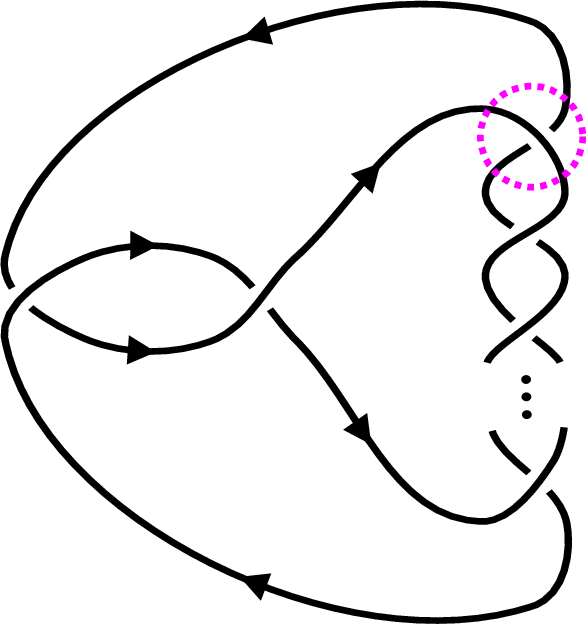}
\caption{$L_+$}
\label{Lplus}
\end{subfigure}
\begin{subfigure}{.3\textwidth}
\centering		\includegraphics[width=0.7\linewidth]{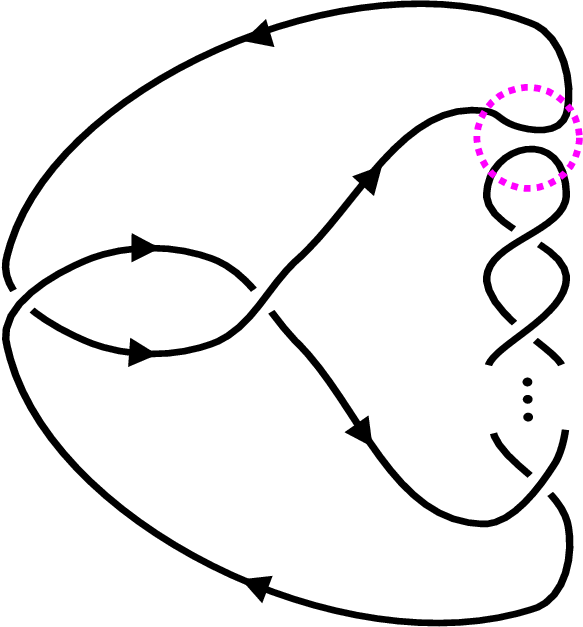}
\caption{$L_0$}
\label{Lzero}
\end{subfigure}
\caption{Skein diagrams of \( P(1,1,2(k+1)+1) \).}
\label{Skeinrelationproof1_1}
\end{figure}

Similarly it is proven for \( n \) an odd negative integer. 

\

Now, suppose that \( n \) is a positive even integer, then we can write \( n = 2q \), where \( q \in \mathbb{Z}^+ \). Let us apply induction on \( q \). For \( q = 1 \), we have \( n = 2 \), which corresponds to the knot \( P(1,1,2) \). Therefore, we have
\[
    \nabla_{P(1,1,2)}(z) = 1 - z^2 = 1 - \left( \frac{2}{2} \right) z^2.
\]
Assume the statement is true for \( q = k \), i.e., for \( P(1,1,2k) \). Consider \( q = k + 1 \). In this case, \( n = 2(k+1) \), and we will prove the formula for \( P(1,1,2(k+1)) \). Let us take a crossing of the knot \( P(1,1,2(k+1)) \); we will select the first crossing of the third tangle and determine the diagrams \( L_+ \), \( L_- \), and \( L_0 \). Notice that \( L_- \) is the knot \( P(1,1,2(k+1)) \). In \( L_+ \), due to the crossing change, two half twists are undone in the third tangle, which means that \( L_+ \) is equivalent to the knot \( P(1,1,2k) \), and by the induction hypothesis, we have \( \nabla_{P(1,1,2k)} = 1 - \left( \frac{2k}{2} \right) z^2 \). Finally, \( L_0 \) is the right-handed Hopf link (and recall that $\nabla_{L_0}(z)=z$); see Figure \ref{Skeinrelationproof1_3}. Using the skein relation, we find that the polynomial is:

\begin{equation*}
    \begin{split}
        \nabla_{P(1,1,n)}(z)=\nabla_{P(1,1,2(k+1))}(z)&=\nabla_{L_-}(z)\\&=\nabla_{L_+}(z)-z\nabla_{L_0}(z)\\
        &=\nabla_{P(1,1,2k)}(z)-z\nabla_{L_0}(z)\\
        &=1-\left(\frac{2k}{2}\right)z^2-z(z)\\
        &=1-\left(\frac{2k}{2}\right)z^2-z^2\\
        &=1-\left(\frac{2k}{2}+1\right)z^2\\
        &=1-\left(\frac{2(k+1)}{2}\right)z^2\\
        &=1-\left(\frac{n}{2}\right)z^2.\\
    \end{split}
\end{equation*}

\begin{figure}[ht]
\centering
\begin{subfigure}{.3\textwidth}
\centering		\includegraphics[width=0.7\linewidth]{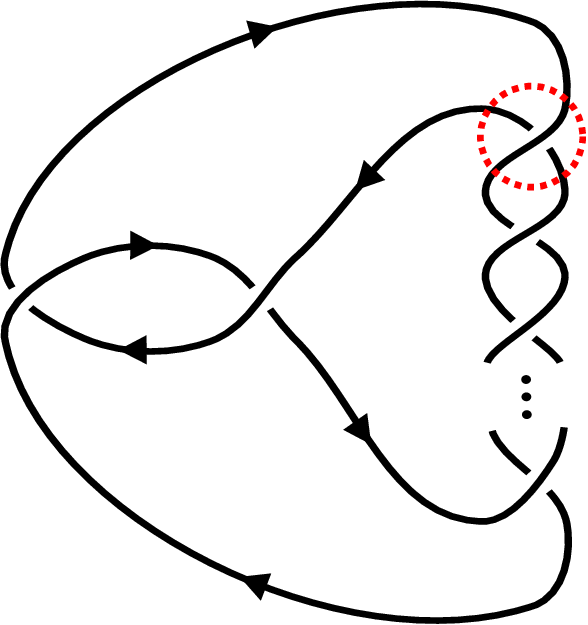}
\caption{$L_-$}
\end{subfigure}
\begin{subfigure}{.3\textwidth}
\centering
\includegraphics[width=0.7\linewidth]{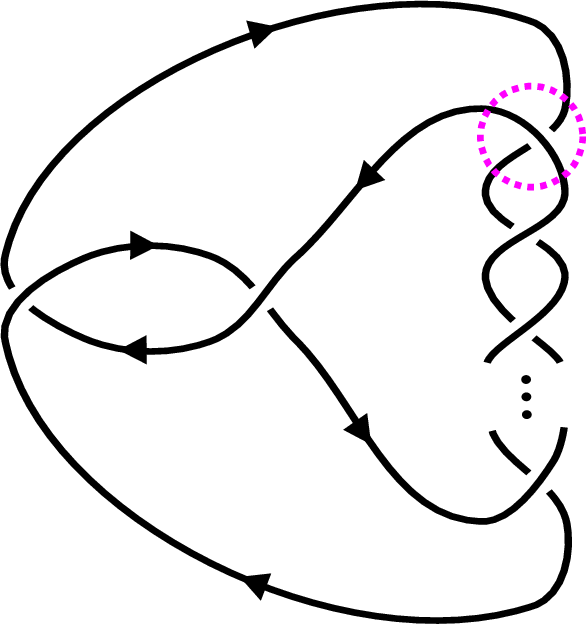}
\caption{$L_+$}
\end{subfigure}
\begin{subfigure}{.3\textwidth}
\centering		\includegraphics[width=0.7\linewidth]{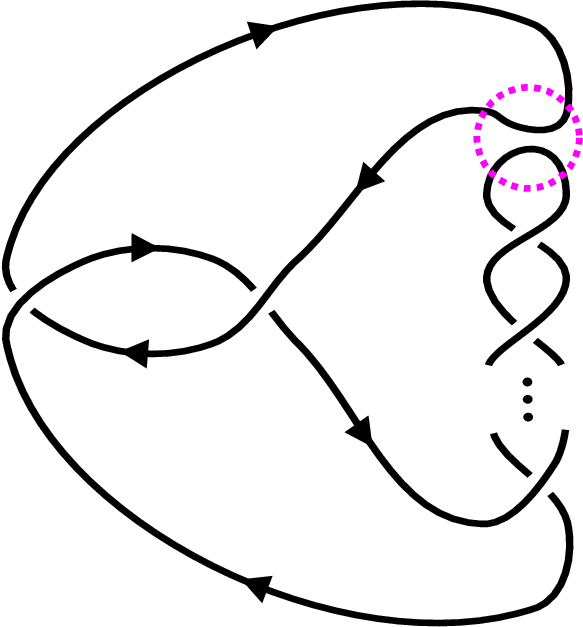}
\caption{$L_0$}
\end{subfigure}
\caption{Skein diagrams of the knot \( P(1,1,2(k+1)) \).}
\label{Skeinrelationproof1_3}
\end{figure}

Similarly it is proven for \( n \) an even negative integer. 
\end{proof}

\section{Kauffman Bracket Polynomial}\label{KBPoly}

In 1985, Louis H. Kauffman introduced a new invariant for unoriented framed links, now known as the Kauffman polynomial in two variables. He initially identified the bracket state summation as a special case of his original two-variable polynomial. At first, Kauffman believed the bracket was a completely new invariant. However, he then realized that the bracket provided a novel and simplified model for the Jones polynomial. For more details on the origins of this polynomial, refer to Chapter 5 in \cite{PBIMW} and \cite{Kau, MV}.

\begin{definition}\ \label{KBPdefi}
\begin{itemize}
\item [(i)] The (reduced) \textbf{Kauffman bracket polynomial} (KBP) is a function from the set of unoriented link diagrams $\mathcal{D}$ to Laurent polynomials with integer coefficients in the variable $A$, $\left< \ \right>: \mathcal{D} \longrightarrow \ \mathbb{Z}[A^{\pm1}]$. The polynomial is characterized by the rules $<\bigcirc>=1$, $<\bigcirc \sqcup K>=(-A^{2}-A^{-2})<K>$, and the skein relation:
$$\left< \KPB \right>=
A\left< \KPC \right> + A^{-1} \left< \KPD \right>.$$

\item [(ii)]	Let $D$ be an unoriented link diagram and let $cr(D)$ be its crossings set. A \textbf{Kauffman state} s, of $D$, is a function $s: cr(D)\longrightarrow \{A, B\}$. This function is understood as an assignment of a \textbf{marker} to each crossing according to the convention illustrated in Figure \ref{fig:kaufstatefunction}. Denote by \textit{KS} the set of all Kauffman states. Moreover, every marker yields a natural \textbf{smoothing} of the crossing as shown in  Figure \ref{fig:kaufstatefunction}. 
\begin{figure}[ht]
\centering
\includegraphics[width=0.95\linewidth]{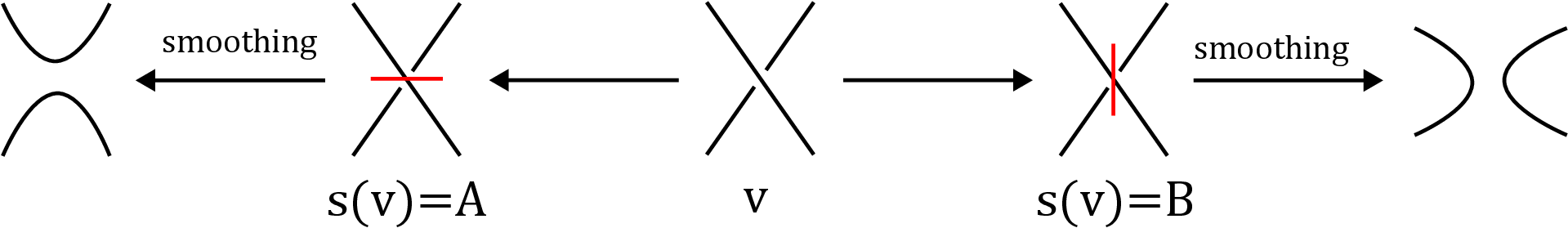}
\caption{Markers at a crossing v of $D$ and their corresponding smoothing.}
\label{fig:kaufstatefunction}
\end{figure}
	
\end{itemize}
	    
\end{definition}

The KBP of a link diagram $D$ is given by the state sum  formula: 
$$<D>=\sum_{s \ \in \ KS}^{}A^{\mid s^{-1}(A) \mid - \mid s^{-1}(B) \mid }(-A^{2}-A^{-2})^{\mid D_{s} \mid -1},$$
where $D_{s}$ denotes the system of circles obtained after smoothing all crossings of $D$ according to the markers of $s$, and $|D_{s}|$ denotes the number of circles in the system.

\begin{example}
We calculate the bracket polynomial of the pretzel knot $P(1,1,1)$. Figure \ref{StatesPretzcel(1,1,1)} shows a diagram of the pretzel knot $P(1,1,1)$ with its crossings labeled $1$, $2$, and $3$. This figure also shows all the states of the knot. The states are labeled with three letters indicating the marker at each crossing, following the order of the labelling. For instance, $ABA$ indicates that crossings $1$ and $3$ are given an $A$ marker while crossing $2$ is given a $B$ marker.

\begin{figure}[H]
    \centering
\includegraphics[width=0.6\textwidth]{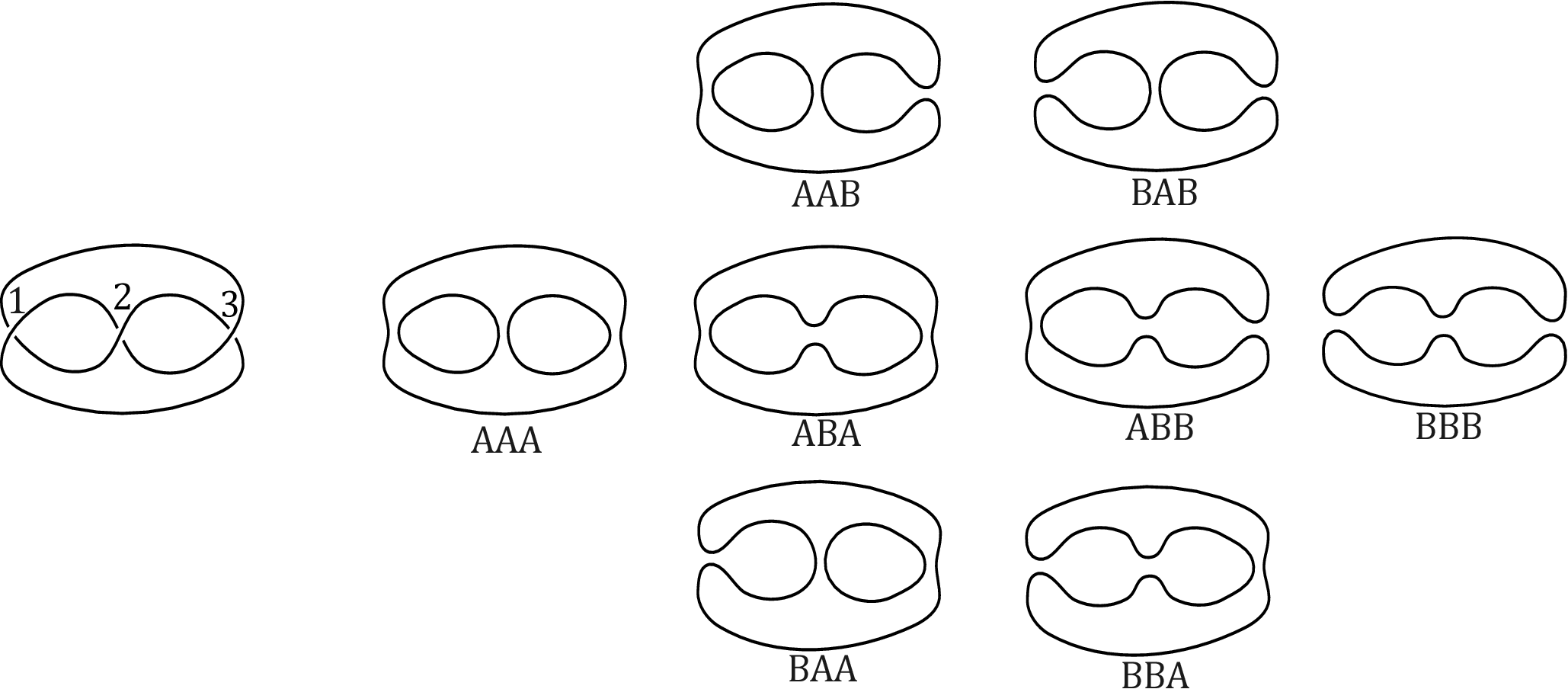}
\caption{States of the pretzel knot $P(1, 1, 1)$.}
\label{StatesPretzcel(1,1,1)}
\end{figure}

Using the KBP formula, we know the contribution from each state:
\begin{gather*}
    AAA \to A^{3-0}(-A^2-A^{-2})^{3-1}=A^3(-A^2-A^{-2})^{2},\\
    AAB, ABA, BAA \to A^{2-1}(-A^2-A^{-2})^{2-1}=A(-A^2-A^{-2}),\\
    BAB, ABB, BBA \to A^{1-2}(-A^2-A^{-2})^{1-1}=A^{-1},\\
    BBB \to A^{0-3}(-A^2-A^{-2})^{2-1}=A^{-3}(-A^2-A^{-2}).
\end{gather*}
Thus the KBP of $P(1, 1, 1)$ is given by: \begin{gather*}
    A^3(-A^2-A^{-2})^{2}+3A(-A^2-A^{-2})+3A^{-1}+A^{-3}(-A^2-A^{-2})=A^7-A^3-A^{-5}.
\end{gather*}    
\end{example}

\subsection{Bracket Polynomial of Pretzel links \texorpdfstring{$\boldsymbol{P(1,1,n)}$}{P(1,1,n)}.}

In this section we present a formula for the calculation of the bracket polynomial of pretzel links $P(1,1,n)$ for a positive integer $n$. First, we explore how to know the number of circles in the states according to the markers. 

\begin{lemma}\label{LemmaForTheoremBracket}
Let $v_1$ and $v_2$ the crossings of the first and second tangle of the pretzel link $P(1,1,n)$, respectively. Let $n>1$ be a positive integer. The number of circles in the state $s$, $|D_s|$, can be determined depending on the number of $B$ markers in the state, as follows.

\begin{itemize}
       \item If $|s^{-1}(B)|=0$, then $|D_s|=3$.
       \item If $|s^{-1}(B)|=1$, then $|D_s|=2$.
       \item If $s^{-1}(v_1)=B$ or $s^{-1}(v_2)=B$ with $|s^{-1}(B)|>1$, then $|D_s|=|s^{-1}(B)|-1$.
       \item If $s^{-1}(v_1) \neq B$ and $s^{-1}(v_2) \neq B$ with $|s^{-1}(B)|>1$, then $|D_s|=|s^{-1}(B)|+1$.
   \end{itemize}
\end{lemma}

\begin{proof} All crossings in the rectangles of the figures in the proof are given $A$ markers.
\begin{itemize}
    \item Suppose $s^{-1}(B)=0$. In this case, all the crossings are given $A$ markers and we obtain three circles in the state, as shown in Figure \ref{AllAmarkersbFirstCase}.

\begin{figure}[ht]
\centering
\begin{subfigure}{.45\textwidth}
\centering		\includegraphics[width=0.25\linewidth]{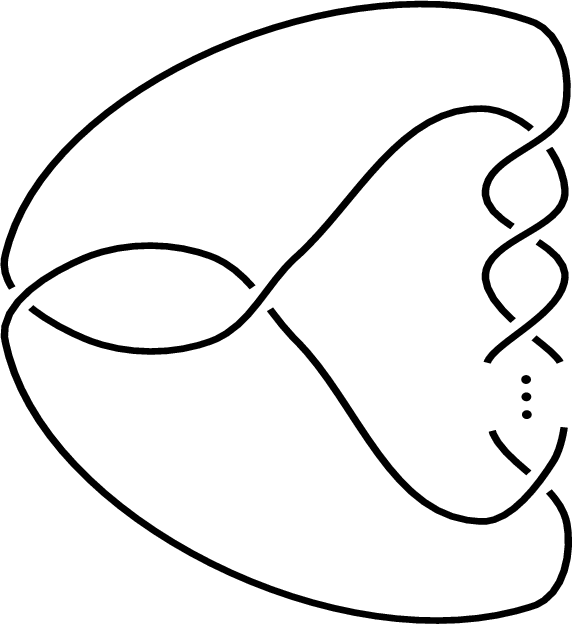}
\caption{$P(1,1,n)$, $n>1$.}
\label{AllAmarkersaFirstCase}
\end{subfigure}
\begin{subfigure}{.45\textwidth}
\centering
\includegraphics[width=0.25\linewidth]{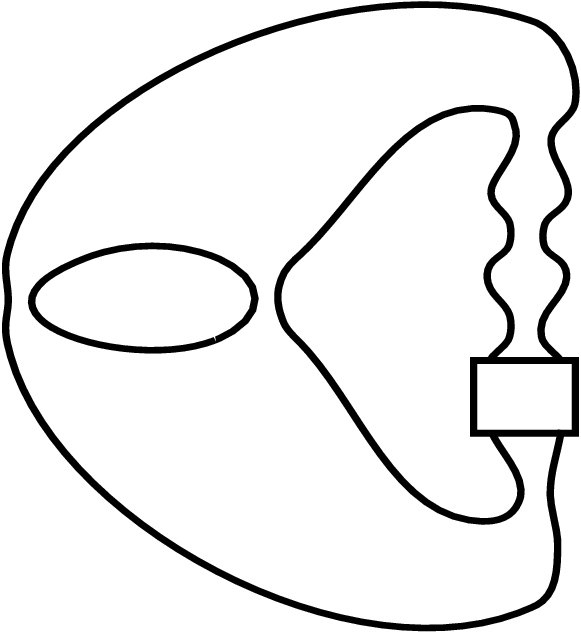}
\caption{State with $s^{-1}(B)=0$.}
\label{AllAmarkersbFirstCase}
\end{subfigure}
\caption{State with all $A$ markers.}
\label{FirstCaseLemma}
\end{figure}

\item Suppose $s^{-1}(B)=1$. In this case, we have three possibilities, depending on which tangle this $B$ marker is located. Observe from Figure \ref{SecondCaseLemma} that $|D_s|=2$.

\begin{figure}[ht]
\centering
\begin{subfigure}{.3\textwidth}
\centering		\includegraphics[width=0.35\linewidth]{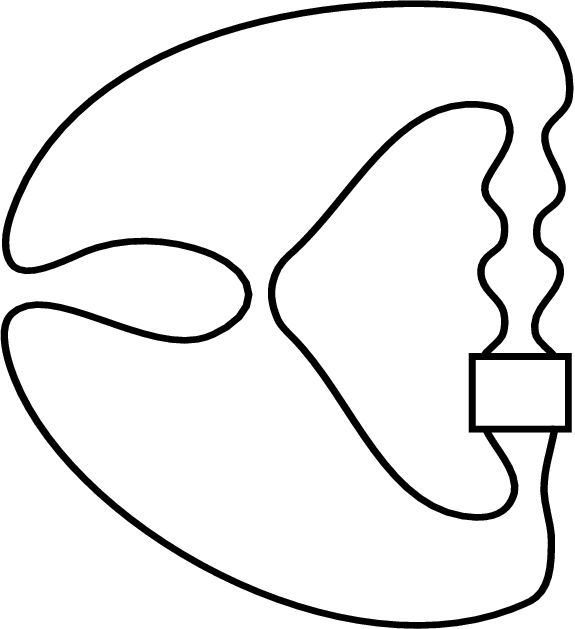}
\caption{$s^{-1}(v_1)=B$.}
\label{AllAmarkersaSecondCase}
\end{subfigure}
\begin{subfigure}{.3\textwidth}
\centering
\includegraphics[width=0.35\linewidth]{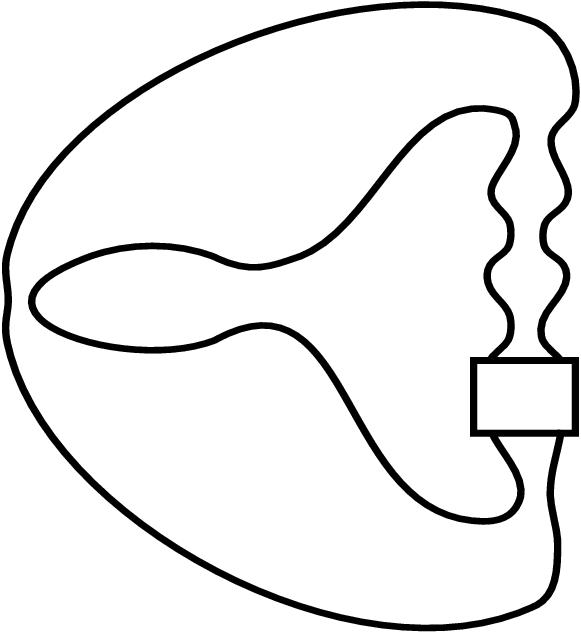}
\caption{$s^{-1}(v_2)=B$.}
\label{AllAmarkersbSecondCase}
\end{subfigure}
\begin{subfigure}{.3\textwidth}
\centering		\includegraphics[width=0.35\linewidth]{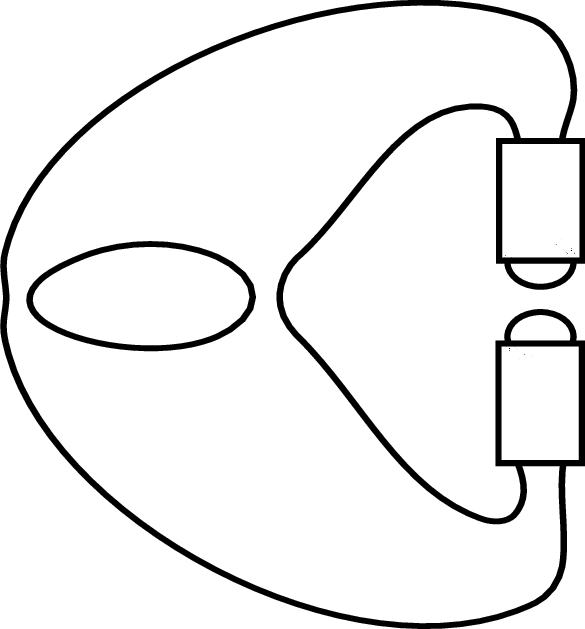}
\caption{$B$ marker in the third tangle.}
\label{BMarkerThirdTangle}
\end{subfigure}
\caption{$s^{-1}(B)=1$.}
\label{SecondCaseLemma}
\end{figure}

\item Suppose $s^{-1}(v_1)=B$ or $s^{-1}(v_2)=B$ with $|s^{-1}(B)|>1$. Observe that the crossings of the first two tangles are defined in the sense that either one has a $B$ marker or both are given $B$ markers. We study the situation in the third tangle with $m$ crossings given a $B$ marker and show that in this case $m-1$ circles are obtained. Let us use an argument by induction. Take as the initial case $m=1$. Observe that no circles (or $m-1$ circles) are obtained in the third column; see Figure \ref{ThirdTangle}.

\begin{figure}[ht]
\centering
\begin{subfigure}{.45\textwidth}
\centering		\includegraphics[width=0.3\linewidth]{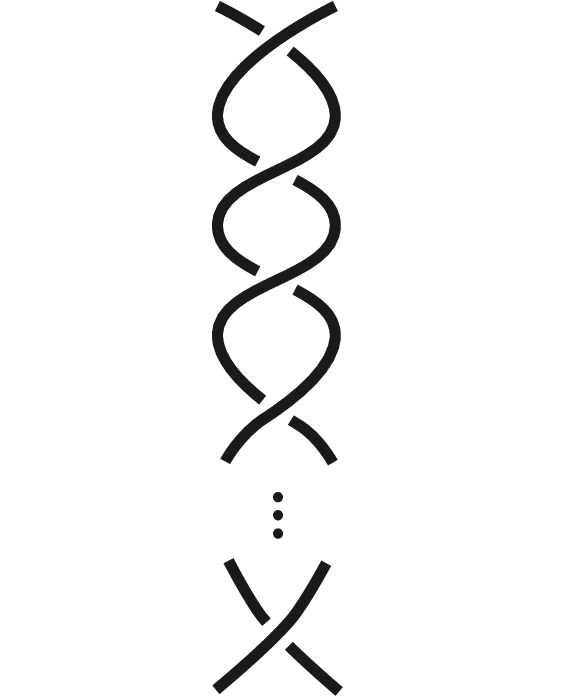}
\caption{Tangle with positive crossings.}
\label{TangleAllPos}
\end{subfigure}
\begin{subfigure}{.45\textwidth}
\centering
\includegraphics[width=0.2\linewidth]{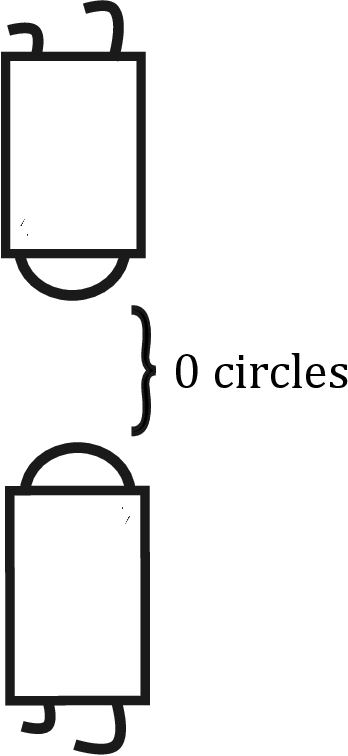}
\caption{No circles if $m=1$.}
\label{ThirdTangle}
\end{subfigure}
\caption{Argument for $m=1$.}
\label{ThirdCaseLemma1}
\end{figure}

As the inductive hypothesis, suppose that for $m=k$ the number of circles obtained is $k-1$. Then, the case $m=k+1$, means that to the system with $k-1$ circles (see Figure \ref{TanglekB}) in the third tangle, a $B$ marker is added, increasing the number of circles by $1$, as it can be seen in Figure \ref{ThirdCaseLemma2}. So, $k=(k+1)-1=m-1$, as we wanted.

\begin{figure}[H]
\centering
\begin{subfigure}{.45\textwidth}
\centering		\includegraphics[width=0.25\linewidth]{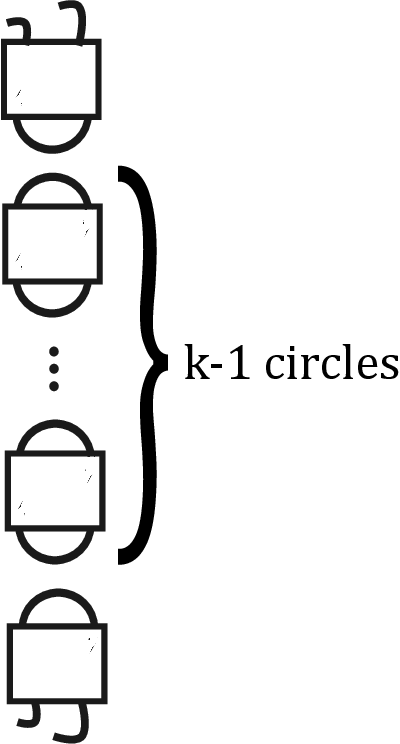}
\caption{Tangle with $k$ crossings given a $B$ marker.}
\label{TanglekB}
\end{subfigure}
\begin{subfigure}{.45\textwidth}
\centering
\includegraphics[width=0.2\linewidth]{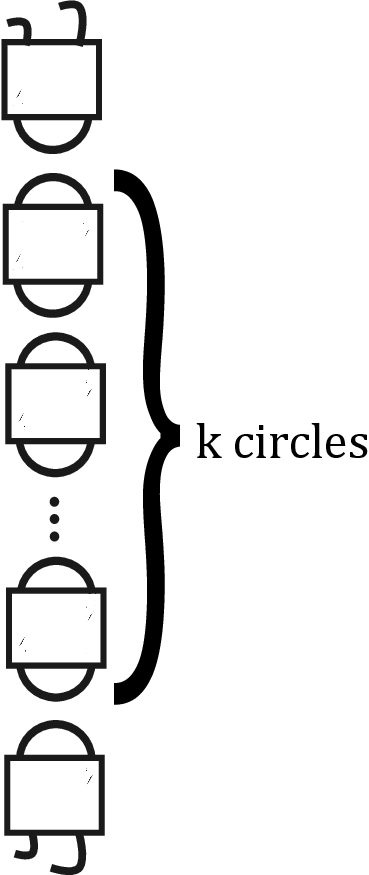}
\caption{Tangle with $k+1$ crossings given a $B$ marker.}
\label{k+1circles}
\end{subfigure}
\caption{Number of circles depending on the number of $B$ markers.}
\label{ThirdCaseLemma2}
\end{figure}

Now we are able to show that when $j=|s^{-1}(B)|>1$, it holds that $|D_s|=|s^{-1}(B)|-1$. Suppose first that $s^{-1}(v_1)=B$ and there are $m$ crossings with a $B$ marker in the third column. Then there are $m-1$ circles in the third tangle and there is another circle formed because of the $B$ marker at $v_1$; see Figure \ref{ThirdCaseONE}. Thus, $|D_s|=m-1+1=(m+1)-1=j-1=|s^{-1}(B)|-1$. Analogously when $s^{-1}(v_2)=B$, $|D_s|=|s^{-1}(B)|-1$ (see Figure \ref{ThirdCaseTwo}). Finally, suppose that both crossings $v_1$ and $v_2$ are smoothed according to a $B$ marker, and there are $m$ labels $B$ in the third column. We have $m-1$ circles from the third tangle and two other circles are formed because of the connection coming from the other two tangles, as it can be observed in Figure \ref{ThirdCaseTHREE}. Thus, $|D_s|=m-1+2=(m+2)-1=j-1=|s^{-1}(B)|-1$.

\begin{figure}[ht]
\centering
\begin{subfigure}{.3\textwidth}
\centering		\includegraphics[width=0.65\linewidth]{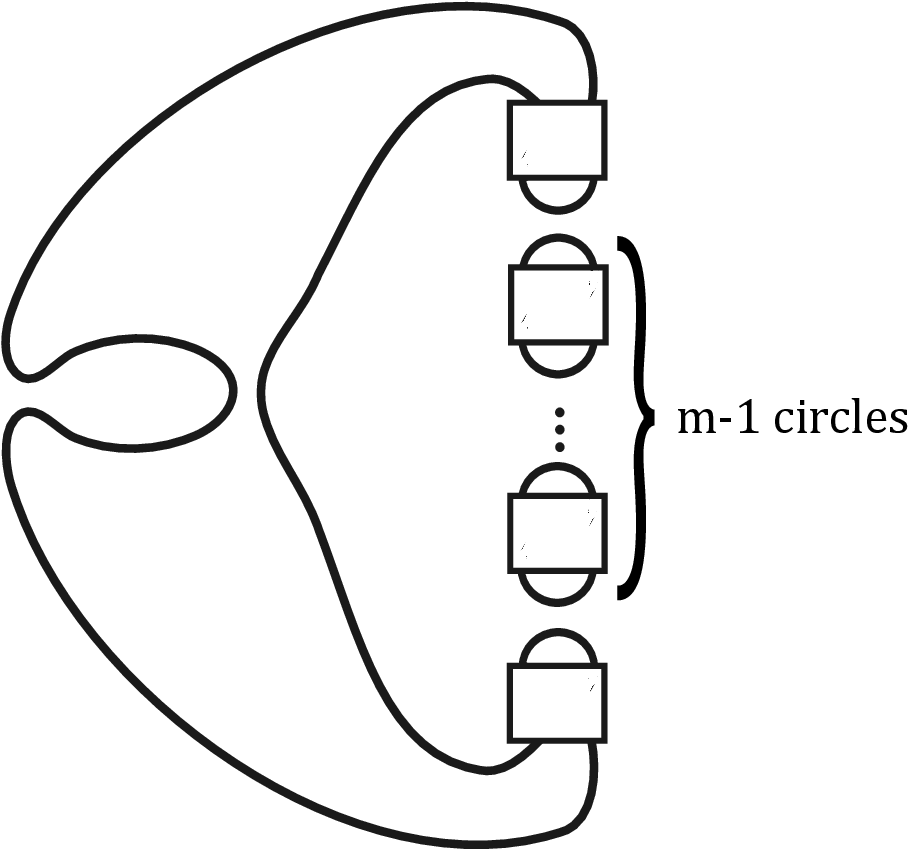}
\caption{$s^{-1}(v_1)=B$. 
}
\label{ThirdCaseONE}
\end{subfigure}
\begin{subfigure}{.3\textwidth}
\centering
\includegraphics[width=0.65\linewidth]{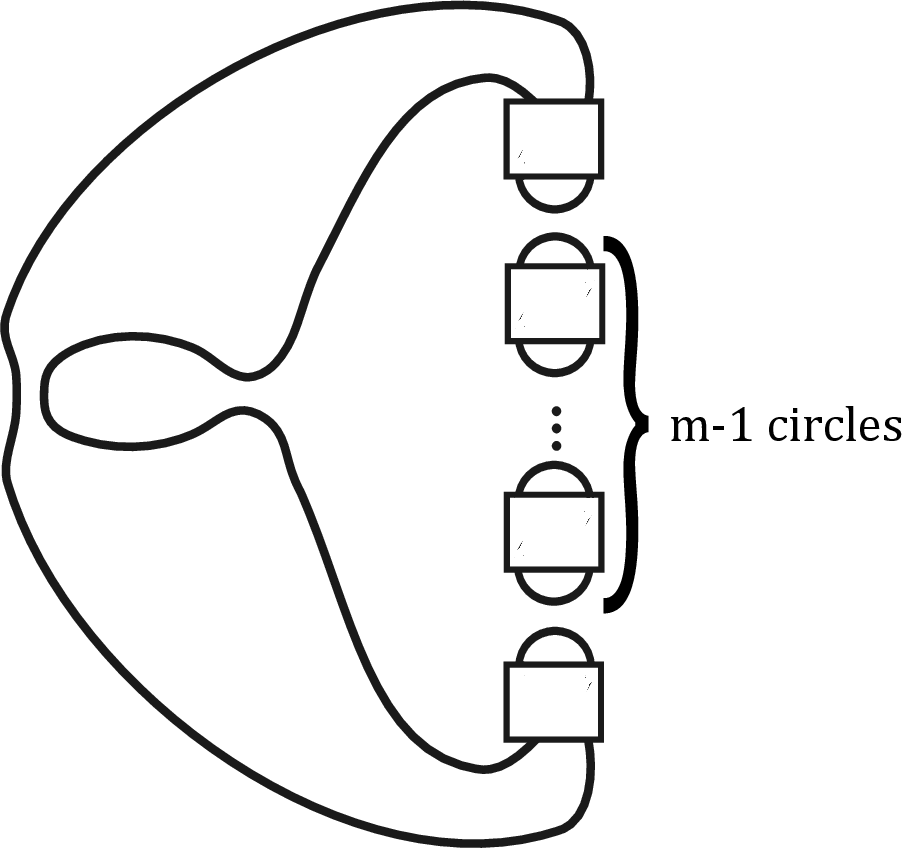}
\caption{$s^{-1}(v_2)=B$. 
}
\label{ThirdCaseTwo}
\end{subfigure}
\begin{subfigure}{.3\textwidth}
\centering		\includegraphics[width=0.65\linewidth]{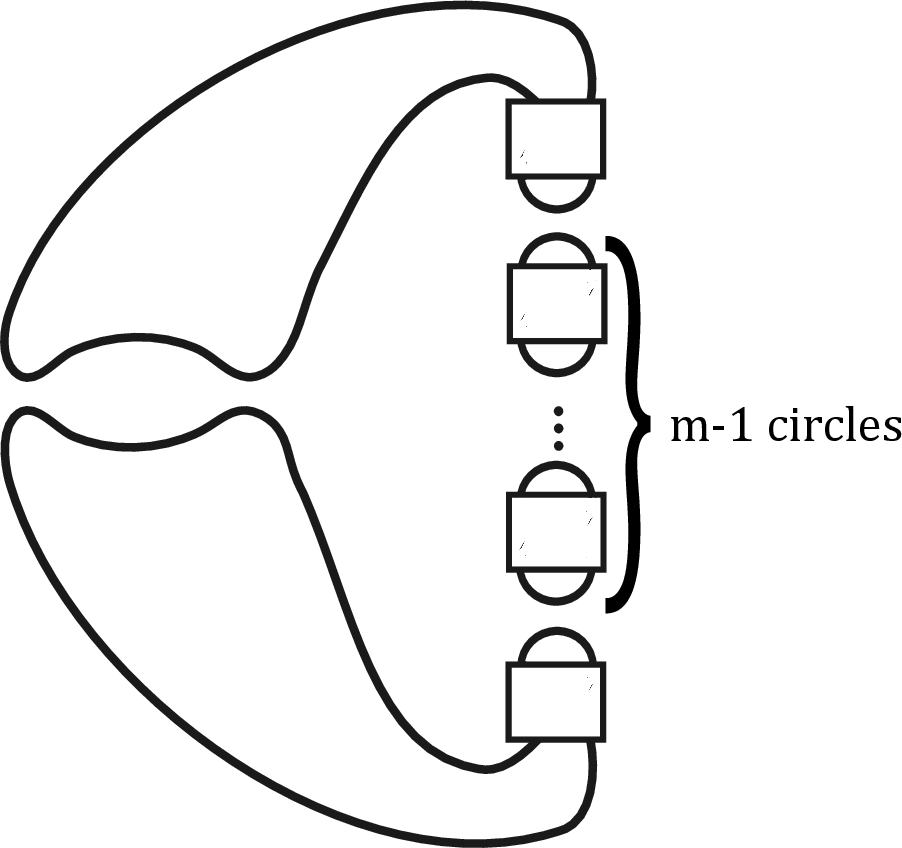}
\caption{$s^{-1}(v_1)=s^{-1}(v_2)=B$.
}
\label{ThirdCaseTHREE}
\end{subfigure}
\caption{Possible combinations of crossings with $B$ markers.}
\label{CombinationsThirdCase}
\end{figure}

\item Suppose that $s^{-1}(v_1) \neq B$ and $s^{-1}(v_2) \neq B$ with $|s^{-1}(B)|>1$. Then all the crossings with $B$ markers are located in the third column. We know there are $m-1$ circles from the third column and, as it can be observed from Figure \ref{LastCase}, there are two other circles generated from the first two tangles. Thus, $|D_s|=m-1+2=m+(2-1)=j+1=|s^{-1}(B)|+1$, as desired.
\begin{figure}[ht]
    \centering
\includegraphics[width=0.2\linewidth]{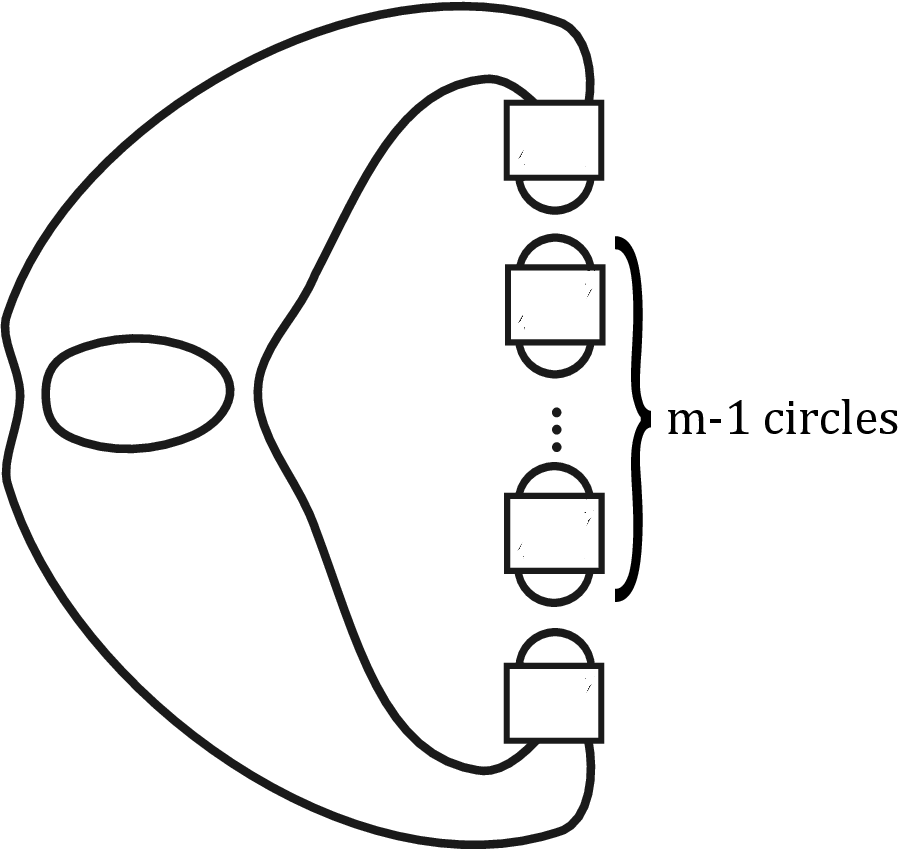}
    \caption{$s^{-1}(v_1) \neq B$ and $s^{-1}(v_2) \neq B$ and there are $m$ crossings with $B$ marker in the third column.}
    \label{LastCase}
\end{figure}    
\end{itemize}
\end{proof}

Now we can state the main result of this section.

\begin{theorem}
    The Kauffman bracket polynomial of a pretzel link $P(1,1,n)$, for $n \in \mathbb{Z}^{+} \setminus {1}$ can be calculated using the following formula:
\begin{equation*}
\begin{split}
 \langle P(1,1,n) \rangle=  & A^{n+2}(-A^2-A^{-2})^2+(n+2)A^n(-A^2-A^{-2}) +  \sum_{j=2}^{n+2} \left[\begin{pmatrix}
        n+1\\
        j-1 
       \end{pmatrix}
       +
       \begin{pmatrix}
        n\\
        j-1 
       \end{pmatrix}\right]
       A^{n+2(1-j)}(-A^2-A^{-2})^{j-2}  \\
 & \hspace{6.7 cm} + \sum_{j=2}^{n} \begin{pmatrix}
        n\\
        j 
       \end{pmatrix}
       A^{n+2(1-j)}(-A^2-A^{-2})^{j}.
\end{split}
\end{equation*}
    
\end{theorem}

\begin{proof}
    Let $K=P(1,1,n)$. We group the states respect to the number of $B$ markers according to Lemma \ref{LemmaForTheoremBracket}, as follows:
\begin{gather*}
    \langle K \rangle=\sum_{s \in KS} A^{|s^{-1}(A)| - |s^{-1}(B)|} (-A^2-A^{-2})^{|D_s|-1}\\
    \\
    \hspace{0.8cm}=\sum_{s \in KS:|s^{-1}(B)|=0} A^{|s^{-1}(A)|}(-A^2-A^{-2})^{|D_s|-1}\\
        \\
        + \sum_{s \in KS:|s^{-1}(B)|=1} A^{|s^{-1}(A)|-1}(-A^2-A^{-2})^{|D_s|-1}\\
        \\
        + \sum_{s \in KS:s^{-1}(v_1) = B \hspace{0.1cm}\lor\hspace{0.1cm} s^{-1}(v_2)=B \hspace{0.1cm} \land \hspace{0.1cm}|s^{-1}(B)|>1} A^{|s^{-1}(A)|-|s^{-1}(B)|}(-A^2-A^{-2})^{|D_s|-1}\\
        \\
        + \sum_{s \in KS:s^{-1}(v_1) \neq B \hspace{0.1cm}\land \hspace{0.1cm} s^{-1}(v_2) \neq B \hspace{0.1cm} \land \hspace{0.1cm}|s^{-1}(B)|>1} A^{|s^{-1}(A)|-|s^{-1}(B)|}(-A^2-A^{-2})^{|D_s|-1}.
\end{gather*}

There is only one state for which $|s^{-1}(B)|=0$ and it follows from Lemma \ref{LemmaForTheoremBracket} that $|D_s|=3$. Moreover, $|s^{-1}(A)|=n+2$ and thus, the first sum is given by $ A^{n+2}(-A^2-A^{-2})^2$.

\ 

The number of states for which $|s^{-1}(B)|=1$, is the number of combinations of having a $B$ marker among all the $n+2$ crossings $\binom{n+2}{1}=n+2$. In this case, $|D_s|=2$ and $|s^{-1}(A)|=n+1$, which implies that $|s^{-1}(A)|-1=n$. Then, the second sum is given by $(n+2)A^n(-A^2-A^{-2})$.

\ 

For the remaining two cases, let $j=|s^{-1}(B)|$. Then we have
\begin{equation*}
    \begin{split}
        |s^{-1}(A)|-|s^{-1}(B)|&=(n+2-j)-j\\
        &=n+2-j-j\\
        &=n+2(1-j).
    \end{split}
\end{equation*}

Consider first the case when $s^{-1}(v_1)=B$. So we fix the first crossing with a $B$ marker and take the combination of having a $B$ marker among the $n+1$ crossings $\binom{n+1}{j-1}$. On the other hand, if $s^{-1}(v_2)=B$ y $s^{-1}(v_1)\neq B$, it means the second crossing is fixed with a $B$ marker and the combination of having a $B$ marker among the crossings on the third column is given by $\binom{n}{j-1}$.

\ 

Observe that the maximum value of j is $n+2$, which means that in the combination we get the case when $j-1>n$. We discard this case (or ay the combination is zero). Now, from Lemma \ref{LemmaForTheoremBracket} we know the number of circles is $|D_s| = |s^{-1}(B)| - 1$, hence the third sum can be written as:

\begin{equation*}
    \sum_{j=2}^{n+2} \left[\begin{pmatrix}
        n+1\\
        j-1 
       \end{pmatrix}
       +
       \begin{pmatrix}
        n\\
        j-1 
       \end{pmatrix}\right]
       A^{n+2(1-j)}(-A^2-A^{-2})^{j-2}.
\end{equation*}

We still need to know the number of states for which $s^{-1}(v_1) \neq B$ y $s^{-1}(v_2) \neq B$ y $|s^{-1}(B)|>1$. The combination of $B$ markers in the third tangle is given by $\binom{n}{j}$ and we know the number of circles in this case is $|D_s|=|s^{-1}(B)| + 1$. Then the fourth sum can be written as:

\begin{equation*}
    \sum_{j=2}^{n} \begin{pmatrix}
        n\\
        j 
       \end{pmatrix}
       A^{n+2(1-j)}(-A^2-A^{-2})^{j}.
\end{equation*}

\end{proof}

The previous result can be extended to the family of knots \( P(1,1,\dots,1,n) \) as follows.

\begin{lemma} \label{LemmaGeneral} 
Consider the pretzel link \( P(1, 1, \ldots, 1, n) \), where \( n \in \mathbb{Z}^+ \) and the number of tangles with a single crossing is \( m \). Let \( p \) be the number of crossings labeled \( A \) in those \( m \) tangles, and \( q \) the number of crossings labeled \( B \) in the tangle with \( n \) crossings. Then, the number of circles in the system \( |D_s| \) is:  
   \begin{itemize}  
       \item $|D_s|=p+q-1$, if $p>0$ and $q>0$.  
       \item $|D_s|=p+1$, if $p >0$ and $q=0$.  
       \item $|D_s|=q+1$, if $p =0$ and $q>0$.  
       \item $|D_s|=1$, if $p =0$ and $q=0$.  
   \end{itemize}  
   \label{Pro_6.14}  
\end{lemma}  

\begin{proof}
First we prove by induction that, given a tangle with positive crossings, $q$ of which are labeled \( B \), exactly \( q - 1 \) circles are generated.

\ 

First let \( q = 1 \) and observe that in this case the number of circles is given by \( 0 = 1 - 1 = q - 1 \) (Figure \ref{FirstLemmaGral}). Suppose now that for \( q = k \), there are \( k - 1 \) circles. Let $q=k+1$. Then, when we increase the number of $B$ markers by one, a new circle is generated, resulting in a total of \( k \) circles; Figure \ref{SecondLemmaGral}. This is \( k = (k + 1) - 1 = q - 1 \). Thus, the number of circles generated in a tangle with positive crossings, of which \( q \) are labeled \( B \), is exactly \( q - 1 \). Similarly, it can be shown that, for a tangle with negative crossings, of which \( p \) are labeled \( A \), exactly \( p - 1 \) circles are generated.

\begin{figure}[ht]
\centering
\begin{subfigure}{.45\textwidth}
\centering		\includegraphics[width=0.3\linewidth]{enredo_1.eps}
\caption{Tangle with $n$ positive crossings.}
\label{TanglePositiveCrossings}
\end{subfigure}
\begin{subfigure}{.45\textwidth}
\centering
\includegraphics[width=0.2\linewidth]{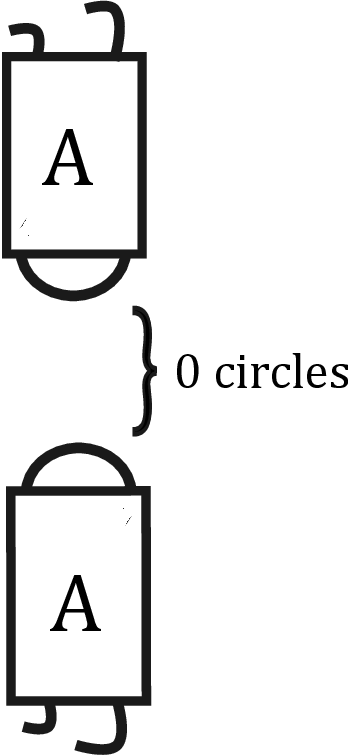}
\caption{Tangle with a crossing with a $B$ marker.}
\label{TanglewithB}
\end{subfigure}
\caption{If $q=1$ then we have $q-1$ circles.}
\label{FirstLemmaGral}
\end{figure}

\begin{figure}[ht]
\centering
\begin{subfigure}{.45\textwidth}
\centering		\includegraphics[width=0.25\linewidth]{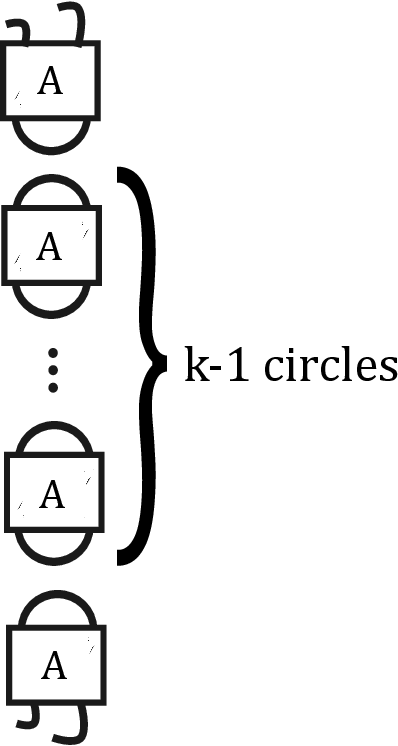}
\caption{Tangle with $k$ crossings with $B$ marker.}
\label{TangleWithkPositCrossings}
\end{subfigure}
\begin{subfigure}{.45\textwidth}
\centering
\includegraphics[width=0.2\linewidth]{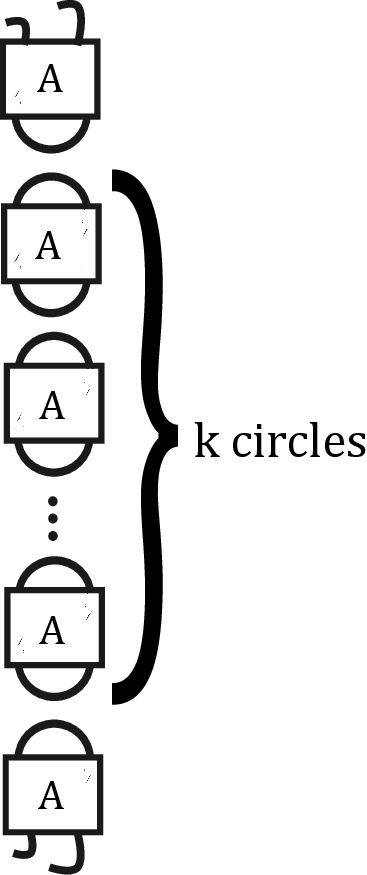}
\caption{Tangle with $k+1$ crossings with $B$ marker.}
\label{TangleWithk+1B}
\end{subfigure}
\caption{Number of circles depending on the number of $B$ markers.}
\label{SecondLemmaGral}
\end{figure}

Now we consider the case when  $p>0$ and $q>0$ we have that there are $q-1$ circles ``inside'' the last tangle and $p-1$ circles ``inside'' the $m$ tangles consisting of one crossing. Moreover, an additional circle is generated given the conecction between the tangles; see Figure \ref{BracketLemma3}. Hence, the number of circles is given by $D_s|=(p-1)+(q-1)+1=p+q-1$.


\begin{figure}[ht]
\centering
\begin{subfigure}{.45\textwidth}
\centering		\includegraphics[width=0.55\linewidth]{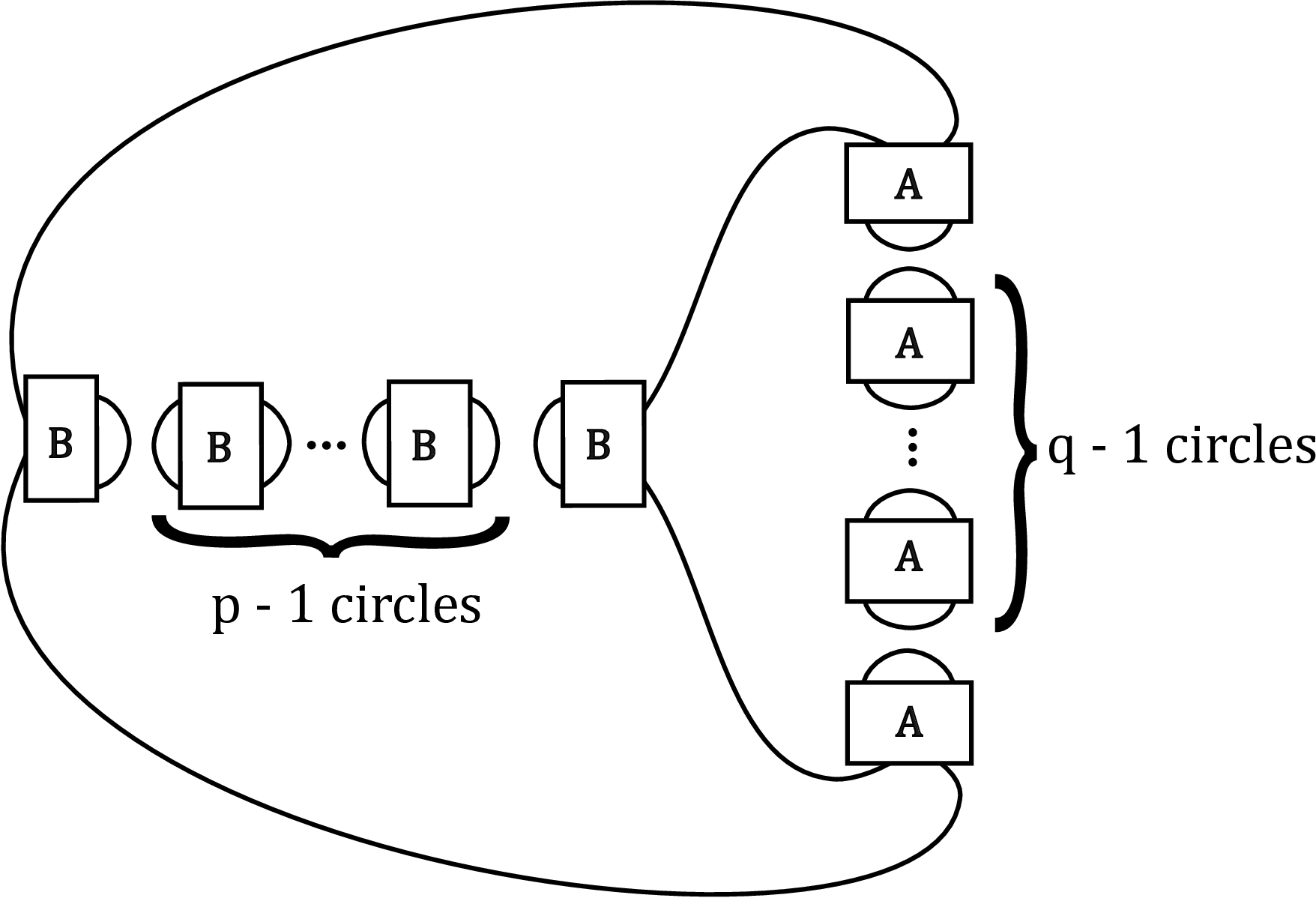}
\caption{$P(1,1,...,1,n)$ with $p>0$ and $q>0$.}
\label{BracketLemma3}
\end{subfigure}
\begin{subfigure}{.45\textwidth}
\centering
\includegraphics[width=0.45\linewidth]{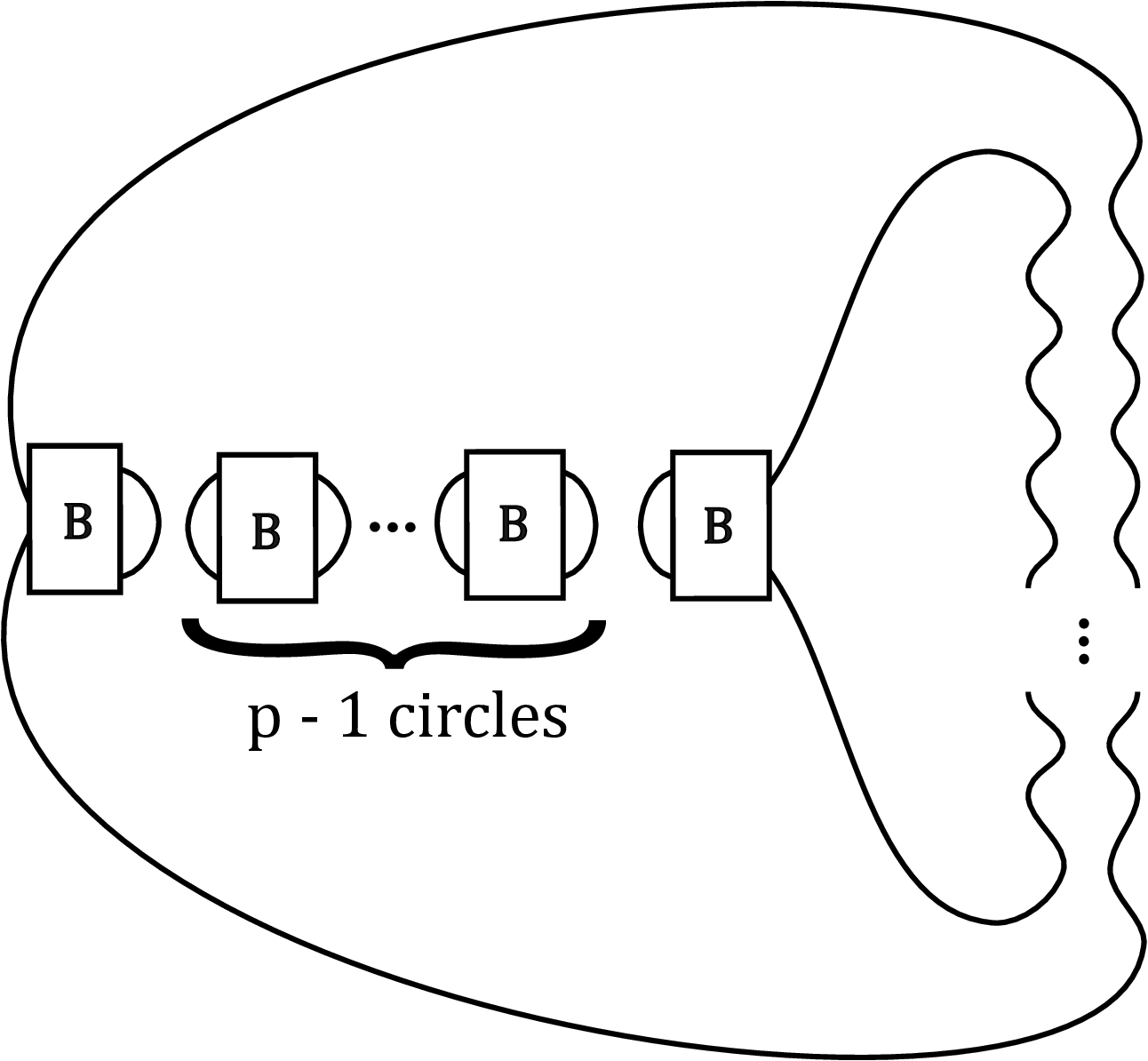}
\caption{$P(1,1,...,1,n)$ with $p>0$ and $q=0$.}
\label{BracketLemma4}
\end{subfigure}
\caption{Number of circles according to the markers.}
\label{ThirdLemmaGral}
\end{figure}

When $p >0$ and $q=0$ we have $p-1$ circles ``inside'' the $m$ tangles with one crossing. Moreover, two circles are generated given the connection between the tangles; see Figure \ref{BracketLemma4}. Hence the number of circles is given by $|D_s|=(p-1)+2=p+1$.

\ 

When $p=0$ and $q>0$ we have $q-1$ circles ``inside'' the last tangle. Moreover, two circles are generated when connecting the tangles; see Figure \ref{BracketLemma5}. Hence, the number of circles is given by $|D_s|=(q-1)+2=q+1$.

\

When $p=0$ and $q=0$ we have that $|D_s|=1$; see Figure \ref{BracketLemma6}.

\begin{figure}[ht]
\centering
\begin{subfigure}{.45\textwidth}
\centering		\includegraphics[width=0.55\linewidth]{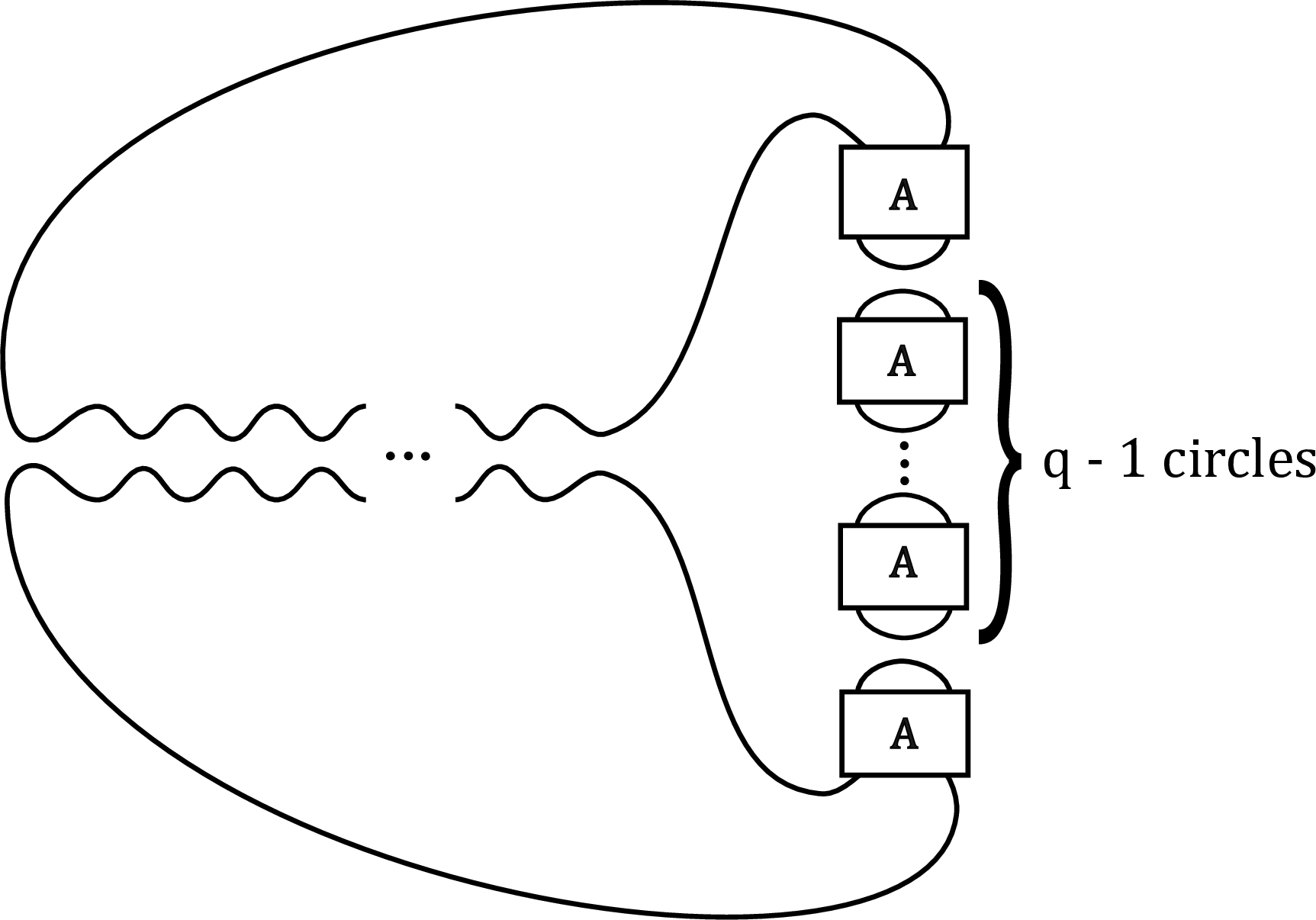}
\caption{$P(1,1,...,1,n)$ with $p=0$ and $q>0$.}
\label{BracketLemma5}
\end{subfigure}
\begin{subfigure}{.45\textwidth}
\centering
\includegraphics[width=0.45\linewidth]{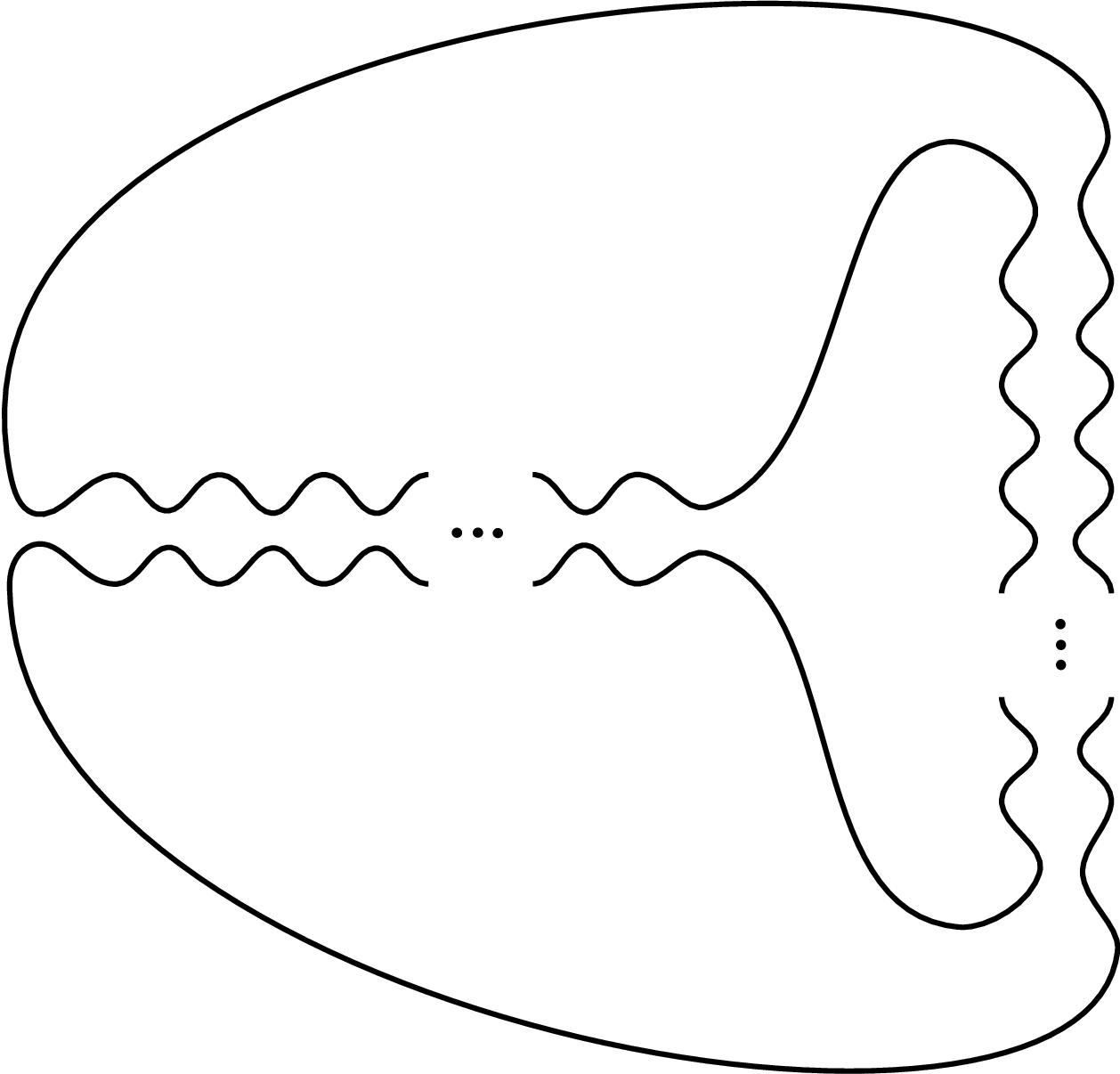}
\caption{$P(1,1,...,1,n)$ with $p=0$ and $q=0$.}
\label{BracketLemma6}
\end{subfigure}
\caption{Number of circles according to the markers.}
\label{FourthLemmaGral}
\end{figure}
\end{proof}

With this lemma, we can prove the following result.

\begin{theorem}  
The bracket polynomial of a pretzel link \( P(1,1,\ldots,1,n) \), with \( n \in \mathbb{Z}^{+} \), where the number of tangles with a single crossing is \( m \), is given by:  
\begin{equation*}  
    \begin{split}  
        \langle P(1,1,\ldots,1,n) \rangle= \sum_{p=1}^{m} \sum_{q=1}^{n} \binom{m}{p} \cdot \binom{n}{q} \left[A^{2(p-q)+n-m}(-A^2-A^{-2})^{p+q-2}\right] \\  
        \\  
        +\sum_{p=1}^{m} \binom{m}{p} \left[A^{2p+n-m}(-A^2-A^{-2})^{p}\right]  
        + \sum_{q=1}^{n} \binom{n}{q} \left[A^{-2q+n-m}(-A^2-A^{-2})^{q}\right] + A^{n-m}.  
    \end{split}  
\end{equation*}  
\end{theorem}

\begin{proof}
We group the states of the pretzel $P(1, 1,...,1, n)$ according to Lemma \ref{LemmaGeneral}. The number of states for each $p>0$ and $q>0$ is:

$$ \binom{m}{p} \cdot \binom{n}{q}.$$

Now, observe that $|s^{-1}(A)|=p+(n-q)$ y $|s^{-1}(B)|=q+(m-p)$. Then $|s^{-1}(A)|-|s^{-1}(B)|=2(p-q)+n-m$. Thus, the first term is obtained:

$$\sum_{p=1}^{m} \sum_{q=1}^{n} \binom{m}{p} \cdot \binom{n}{q}A^{2(p-q)+n-m}(-A^2-A^{-2})^{p+q-2}.$$

Similarly, the number of states for each $p>0$ and $q=0$ is:

$$ \binom{m}{p} \cdot \binom{n}{0}=\binom{m}{p} $$

and the number of states for $q>0$ and $p=0$ is:

$$  \binom{m}{0} \cdot \binom{n}{q}=\binom{n}{q},$$

In this way, the following terms are obtained:

$$\sum_{p=1}^{m} \binom{m}{p}\left[A^{2p+n-m}(-A^2-A^{-2})^{p}\right] + \sum_{q=1}^{n} \binom{n}{q}\left[A^{-2q+n-m}(-A^2-A^{-2})^{q}\right]$$

Finally, we consider the case $p=0$ and $q=0$. Here, there is only one state which gives the term $A^{n-m}$.

\end{proof}

\section{Future Directions}\label{Conclusions}
The study of polynomials in knot theory remains a significant area of interest in the community. They not only serve as an accessible introduction to knot theory research but also offer a wide range of possibilities, from patterns in the calculations to categorifications via homology theories. Exploring other polynomials, such as the HOMFLYPT and plucking polynomials \cite{ILMP}, seems to be an adequate research direction. Furthermore, focusing on specific families of links could enhance our understanding of these polynomials and the properties of those links.


\section*{Acknowledgments}
This article partially results from work done towards the degree project of the first author at the Universidad Nacional Aut\'onoma de Honduras (UNAH). The second author acknowledges the support of the National Science Foundation through Grant DMS-2212736.


\end{document}